\documentclass[a4paper, 10pt,oneside]{amsart}
\usepackage{amsmath, amssymb, amsfonts}
\usepackage[all]{xy}
\newtheorem{tm}{Theorem}
\newtheorem{lm}[tm]{Lemma}
\newtheorem{prop}[tm]{Proposition}
\newtheorem{kor}[tm]{Corollary}

\newenvironment{dokaz}
{\noindent\emph{Proof:}\ }
{\hfill $\blacksquare$}

\newcommand{\Z}
{{\mathbb Z}}
\newcommand{\N}
{{\mathbb N}}
\newcommand{\C}
{{\mathbb C}}

\newcommand{\g}
{{\mathfrak g}}
\newcommand{\lgg}
{{\mathfrak l}}

\newcommand{\gt}
{\tilde{{\mathfrak g}}}

\newcommand{\hz}
{\hat{\mathfrak h}_\Z}
\newcommand{\he}
{{\mathfrak h}^e}
\newcommand{\nt}
{\tilde{{\mathfrak n}}}
\newcommand{\h}
{{\mathfrak h}}
\newcommand{\n}
{{\mathfrak n}}
\newcommand{\gsl}
{{\mathfrak sl}}

\begin{document}

\author{Goran Trup\v{c}evi\'{c}}
\title{Combinatorial bases of Feigin-Stoyanovsky's type subspaces of higher-level standard $\tilde{\gsl}(\ell+1,\C)$-modules}
\address{Department of Mathematics, University of Zagreb, Bijeni\v cka 30, Zagreb, Croatia}
\curraddr{}
\email{gtrup@math.hr}
\thanks{}
\subjclass[2000]{Primary 17B67; Secondary 17B69, 05A19.\\ \indent Partially supported by the Ministry of Science and Technology of the Republic of Croatia, Project ID 037-0372794-2806}
\date{}
\dedicatory{}

\begin{abstract}

Let $\tilde{\mathfrak g}$ be an affine Lie algebra of the type
$A_\ell^{(1)}$.
We find a combinatorial basis of Feigin-Stoyanovsky's type subspace $W(\Lambda)$ given
in terms of difference and initial conditions. Linear
independence of the generating set is proved  inductively by using
coefficients of intertwining operators. A basis of $L(\Lambda)$ is obtained as an ``inductive limit'' of the basis of $W(\Lambda)$.
\end{abstract}

\maketitle


\section{Introduction}

Vertex-operator construction of affine Lie algebras has been used
to prove and to give a representation-theoretic interpretation of
Rogers-Ramanujan-type combinatorial identities. This approach was
initiated by J.Lepowsky and R.Wilson ([LW]), and was continued in
the works of Lepowsky, M.Primc, A.Meurman and others (cf. [LP],
[MP]). An important part of this program  was to find monomial
bases of standard modules for affine Lie algebras, or some of its
subspaces. Knowledge of bases was then used to calculate the
character of these spaces, which gave the sum side in the
Rogers-Ramanujan-type partition identities.

Later,  B.Feigin and A.Stoyanovsky considered what they called a
principal subspace of the basic $\tilde{\gsl}(2,\C)$-module
([FS]). They have constructed a combinatorial basis of this
subspace, which again gave the sum side of a Rogers-Ramanujan
identity. G.Georgiev extended character formulas obtained by
Feigin and Stoyanovsky to a family of standard
$\tilde{\gsl}(\ell+1,\C)$-modules ([G]). He explicitly constructed
combinatorial bases for these subspaces and in the proof of linear
independence, he used intertwining operators from [DL]. More
recently, S.Capparelli, J.Lepowsky and A.Milas laid out a program
to interpret and obtain Rogers-Ramanujan-type recursions in the
setting of vertex-operator algebras and affine Lie algebras. They
used intertwining operators to construct exact sequences between
different principal subspaces for $\tilde{\gsl}(2,\C)$ and in this
way, obtained Rogers-Ramanujan and Rogers-Selberg recursions for
characters of these subspaces ([CLM1,CLM2]). As a continuation of
this program, C.Calinescu obtained Rogers-Ramanujan-type
recursions for some classes of standard modules for
$\tilde{\gsl}(\ell+1,\C)$ ([C1,C2]), and Calinescu, Lepowsky and
Milas provided new proofs of presentation theorems for principal
subspaces for $\tilde{\gsl}(2,\C)$ ([CalLM1,CalLM2]).

In parallel with these developments, M.Primc studied similar
subspaces of standard modules for different affine Lie algebras
([P1,P2]), which he later called Feigin-Stoyanovsky's type
subspaces ([P3]). He used bases of these subspaces to construct
from them bases of the whole standard modules. For
$\tilde{\gsl}(3,\C)$, these bases were parameterized
 by $(k,\ell+1)$-admissible configurations, a combinatorial objects that were introduced and further studied in [FJLMM] and [FJMMT].
In [P3], Primc proved linear independence of the spanning set by
using Capparelli-Lepowsky-Milas' approach via intertwining
operators and a description of the basis from [FJLMM]. These
operators and a description of basis were used by M.Jerkovi\' c to
obtain exact sequences of Feigin-Stoyanovsky's type subspaces and
recurrence relations for the corresponding characters ([J1]). By
solving these relations in the $\tilde{\gsl}(3,\C)$-case he was
able to generalize character formulas from [FJLMM].

In our previous paper ([T]) we have used ideas of Georgiev,
Capparelli, Lepowsky and Milas, and of Primc to construct and
prove linear independence of the spanning set for a
Feigin-Stoyanovsky's type subspace for all basic modules for
$\tilde{\gsl}(\ell+1,\C)$. In this paper we generalize this result
to higher-level standard modules for $\tilde{\gsl}(\ell+1,\C)$.

Let $\g=\gsl(\ell+1,\C)$ be a simple complex Lie algebra of type $A_\ell$, $\h\subset\g$ its Cartan
subalgebra, $R$ the corresponding root system. Then one has a root
decomposition $\g=\h+\sum_{\alpha\in R}\g_\alpha$. Fix root vectors
$x_\alpha\in\g_\alpha$. Let $\langle
\cdot,\cdot\rangle$ be a normalized invariant bilinear form on
$\g$, and by the same symbol denote the induced form on $\g^*$.

\begin{figure}[ht] \caption{The set of colors $\Gamma$} \label{Gamma_fig}
\begin{center}\begin{picture}(245,140)(-25,-10) \thicklines
\put(0,0){\line(1,0){180}} \put(0,0){\line(0,1){120}}
\put(180,120){\line(0,-1){120}} \put(180,120){\line(-1,0){180}}
\put(-6,111){$\scriptstyle 1$}
\put(-6,99){$\scriptstyle 2$} \put(-8,3){$\scriptstyle m$}
\put(3,-8){$\scriptstyle m$} \put(12,-8){$\scriptstyle m
{\scriptscriptstyle +} 1$} \put(174,-8){$\scriptstyle \ell$}

 \linethickness{.075mm} \multiput(0,108)(4,0){45}{\line(1,0){2}}
\multiput(0,96)(4,0){45}{\line(1,0){2}}
\multiput(0,12)(4,0){45}{\line(1,0){2}}
\multiput(12,0)(0,4){30}{\line(0,1){2}}
\multiput(24,0)(0,4){30}{\line(0,1){2}}
\multiput(168,0)(0,4){30}{\line(0,1){2}}
\multiput(0,54)(4,0){45}{\line(1,0){2}}
\multiput(0,66)(4,0){45}{\line(1,0){2}}
\multiput(90,0)(0,4){30}{\line(0,1){2}}
\multiput(102,0)(0,4){30}{\line(0,1){2}}

\put(-6,57){$\scriptstyle i$} \put(93,-8){$\scriptstyle j$}
\put(90,57){$\gamma_{ij}$}
\end{picture}\end{center}
\end{figure}
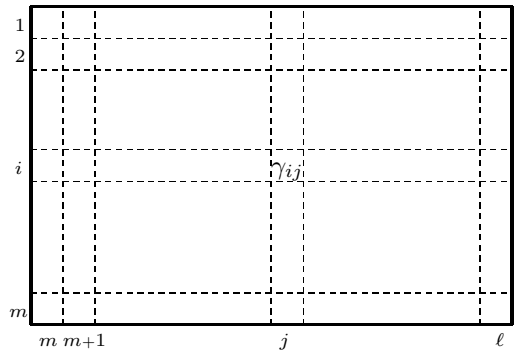

Let $\Pi=\{\alpha_1,\dots,\alpha_\ell\}$ be a basis of the root system $R$, and $\{\omega_1,\dots,\omega_\ell\}$ the
corresponding set of fundamental weights. We identify $\h$ and $\h^*$ in the usual way and fix a fundamental weight
$\omega=\omega_m$. Set $$\Gamma=\{\gamma\in
R\,|\,\langle\gamma,\omega\rangle=1\}=\{\gamma_{ij}\,|\, i=1,\ldots,m; j=m,\ldots,\ell\},$$
where $$ \gamma_{ij}=\alpha_i+\cdots+\alpha_m+\cdots+\alpha_j.$$ Set
 $${\mathfrak g}_{\pm1}  =
\sum_{\alpha \in \pm \Gamma}\, {\mathfrak g}_\alpha,\ {\mathfrak
g}_0 = {\mathfrak h} \oplus\sum_{\langle\alpha,\omega\rangle=0}\,
{\mathfrak g}_\alpha.$$ Then \begin{equation}\label{ZGradG_jed}\g=\g_{-1}\oplus \g_0 \oplus
\g_1\end{equation} is a $\Z$-gradation of
$\g$. The set $\Gamma$ is called {\em the set of colors}. For $\gamma\in \Gamma$, we say that a fixed basis element
$x_\gamma\in\g_\gamma$ {\em is of the color} $\gamma$. The
set of colors $\Gamma$ can be pictured as a rectangle with
row indices $1,\dots,m$ and column indices $m,\dots,\ell$ (see figure \ref{Gamma_fig}).

Affine
Lie algebra associated with $\g$ is $\gt=\g\otimes
\C[t,t^{-1}]\oplus \C c \oplus \C d$, where $c$ is the canonical
central element, and $d$ is the degree operator. Elements
$x_\alpha(n)=x_\alpha\otimes t^n$ are fixed real root vectors.
The $\Z$-gradation \eqref{ZGradG_jed} of $\g$ induces analogous $\Z$-gradation of $\gt$:
$$\gt=\gt_{-1}\oplus \gt_0 \oplus \gt_1,$$
where $\gt_1=\g_1\otimes\C[t,t^{-1}]$ is a commutative Lie
subalgebra with a basis
$$\{x_\gamma(j)\,|\,j\in\Z,\gamma\in\Gamma\}.$$

Let $L(\Lambda)$ be a standard $\gt$-module of level $k=\Lambda(c)$,
with a fixed highest weight vector $v_\Lambda$. A Feigin-Stojanovsky's
type subspace  is a $\gt_1$-submodule of $L(\Lambda)$ generated with
$v_\Lambda$,
$$W(\Lambda)=U(\gt_1)\cdot v_\Lambda\subset L(\Lambda).$$

We find a basis of the Feigin-Stoyanovsky's type subspace
$W(\Lambda)$ consisting of {\em monomial vectors}
$$\{x_{\gamma_1}(-r_1) \cdots x_{\gamma_n}(-r_n) v_\Lambda\,|\,
n\in\Z_+;\gamma_j\in \Gamma, r_j\in\N\}$$ whose {\em monomial parts}
\begin{equation}
\label{MonomPart_def} x_{\gamma_1}(-r_1) \cdots x_{\gamma_t}(-r_n)
\end{equation}
satisfy certain combinatorial conditions, called {\em difference}
and {\em initial conditions}.

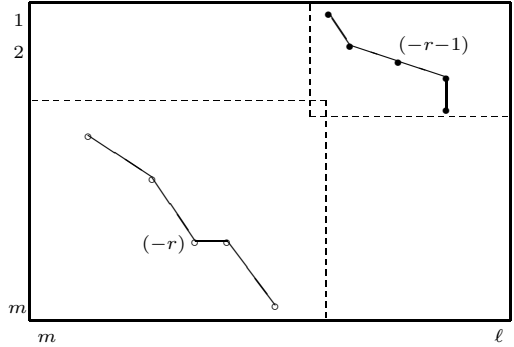
\begin{figure}[ht] \caption{Difference conditions} \label{DCSeq_fig}
\begin{center}\begin{picture}(200,140)(-8,-10) \thicklines
\put(0,0){\line(1,0){180}} \put(0,0){\line(0,1){120}}
\put(180,120){\line(0,-1){120}} \put(180,120){\line(-1,0){180}}
\put(-6,111){$\scriptstyle 1$} \put(-6,99){$\scriptstyle 2$}
\put(-8,3){$\scriptstyle m$} \put(3,-8){$\scriptstyle m$}
\put(174,-8){$\scriptstyle \ell$}

\linethickness{.075mm} \multiput(105,77)(4,0){19}{\line(1,0){2}}
\multiput(105,77)(0,4){11}{\line(0,1){2}}
\multiput(111,83)(-4,0){28}{\line(-1,0){2}}
\multiput(111,83)(0,-4){21}{\line(0,-1){2}}

\thinlines \put(20,68){$\scriptscriptstyle \circ$}
\put(44,52){$\scriptscriptstyle \circ$}
\put(60,28){$\scriptscriptstyle \circ$}
\put(72,28){$\scriptscriptstyle \circ$}
\put(90,4){$\scriptscriptstyle \circ$} \put(22,70){\line(3,-2){24}}
\put(46,54){\line(2,-3){16}} \put(62,30){\line(1,0){12}}
\put(74,30){\line(3,-4){18}}

\put(110,114){$\scriptscriptstyle \bullet$}
\put(118,102){$\scriptscriptstyle \bullet$}
\put(136,96){$\scriptscriptstyle \bullet$}
\put(154,90){$\scriptscriptstyle \bullet$}
\put(154,78){$\scriptscriptstyle \bullet$}
\put(112,116){\line(2,-3){8}} \put(120,104){\line(3,-1){18}}
\put(138,98){\line(3,-1){18}} \put(156,92){\line(0,-1){12}}

\put(138,102){$\scriptstyle (-r-1)$} \put(42,27){$\scriptstyle
(-r)$}

\end{picture}\end{center}
\end{figure}

We'll say that a monomial  \eqref{MonomPart_def} satisfies {\em
difference conditions on} $L(\Lambda)$, if exponents of its factors satisfy the following family of inequalities:
\begin{eqnarray*}
a_{i_1 j_s}^{r+1}+ \dots +a_{i_s j_1}^{r+1}+a_{i_{s+1} j_t}^r+ \dots
+a_{i_t j_{s+1}}^r
\leq k, & & \\
& & \hspace{-58ex} 1\leq i_1 \leq\dots\leq i_s\leq i_{s+1} \leq\dots\leq i_t\leq m, \\
& & \hspace{-58ex}  \ell\geq j_1 \geq\dots\geq j_s \geq j_{s+1}\geq\dots\geq j_t\geq m \\
& & \hspace{-58ex}  (i_\nu,j_\nu)\neq (i_{\nu+1},j_{\nu+1}),\textrm{
for }
\nu=\begin{minipage}[t]{2cm}$1,\dots,s-1,$\\$s+1,\dots,t,$\end{minipage} \nonumber
\end{eqnarray*}
where $a_{ij}^r$ is an exponent of $x_{\gamma_{ij}}(-r)$. We can reformulate this by saying that
for any configuration of colors of elements of degree $-r$ and $-r-1$ of the type pictured on the
figure \ref{DCSeq_fig}, the sum of corresponding exponents must be less than $k+1$.

Similarly, we'll say that a monomial \eqref{MonomPart_def} satisfies {\em
initial conditions on} $L(\Lambda)$ if
\begin{eqnarray*}
& a_{i_1 j_t}^1+ a_{i_2 j_{t-1}}^1 +\dots +a_{i_t j_1}^1\leq
k_0+k_1+\dots+k_{i_t-1}+k_{j_t+1}+\dots+k_\ell, & \\
& & \hspace{-72ex} 1\leq
i_1 \leq i_2\leq\dots\leq i_t\leq m, \\
& &  \hspace{-72ex} \ell\geq j_1 \geq j_2\geq\dots\geq j_t\geq m,\\
& &  \hspace{-72ex} (i_\nu,j_\nu)\neq (i_{\nu+1},j_{\nu+1}).
\nonumber
\end{eqnarray*}

Difference conditions on monomials are obtained by observing
relations between fields $x_\gamma(z),\gamma\in\Gamma$ on
$L(\Lambda)$, while initial conditions are consequences of relations for some
modules of lower level.

In the case $k=1$, these conditions are equivalent to the conditions
we have obtained in [T]. Since a standard module of level $k$ can be found inside a $k$-fold tensor product of modules of level $1$,
it is natural to ask  can a monomial \eqref{MonomPart_def} be factorized in such a way so that each factor satisfies difference and initial
conditions on the corresponding level $1$ module. By combinatorial arguments, we show that the answer to this question is affirmative. This
enables us to use tensor products of coefficients of intertwining operators that were constructed in [T] and to inductively prove the linear
independence of the set of
monomial vectors $x_{\gamma_1}(-r_1) \cdots x_{\gamma_n}(-r_n) v_\Lambda$
whose monomial parts satisfy difference and initial  conditions. Thus, this set is a basis of $W(\Lambda)$.

Following the approach of Primc ([P1],[P2]) we construct the basis of the whole standard module  $L(\Lambda)$ as an ``inductive limit'' of
the basis of $W(\Lambda)$.

\section{Affine Lie algebras}

For $\ell\in \N$, let $$\g=\gsl(\ell+1,\C),$$ a simple Lie algebra
of the type $A_\ell$. Let $\h\subset\g$ be a Cartan subalgebra of
$\g$ and $R$ the corresponding root system. Fix a basis
$\Pi=\{\alpha_1,\dots,\alpha_\ell\}$ of $R$. Then we have the
triangular decomposition $\g=\n_-\oplus \h \oplus \n_+$. By $R_+$
and $R_-$ we denote sets of positive and negative roots, and let
$\theta$ be the maximal root. Let $\langle
x,y\rangle=\textrm{tr\,}xy$ be a normalized invariant bilinear form on
$\g$; via $\langle \cdot,\cdot\rangle$ we have an identification
$\nu:\h\to\h^*$. For each root $\alpha$ fix a root vector $x_\alpha\in\g_\alpha$.

Let $\{\omega_1,\dots,\omega_\ell\}$ be the set of
fundamental weights of $\g$,
$\langle\omega_i,\alpha_j\rangle=\delta_{ij},\,i,j=1,\dots,\ell$. Denote by $Q=\sum_{i=1}^\ell\Z \alpha_i$
the root lattice, and by $P=\sum_{i=1}^\ell\Z \omega_i$ the weight lattice
of $\g$.

Denote by $\gt$ the associated affine Lie algebra
$$\gt=\g\otimes \C[t,t^{-1}]\oplus \C c \oplus \C d.$$
Set $x(j)=x\otimes t^j$ for $x\in\g,j\in\Z$. Commutation relations
are then given by
\begin{eqnarray*}
\ [c,\gt] & = & 0, \\
\ [d,x(j)] & = & j x(j), \\
\ [x(i),y(j)] & = & [x,y](i+j)+ i\langle x,y\rangle \delta_{i+j,0}c.
\end{eqnarray*}


Set $\he=\h\oplus\C c \oplus\C d,\, \nt_\pm=\g\otimes t^{\pm 1}\C
[t^{\pm 1}]\oplus \n_\pm$. Then we also have the triangular
decomposition $\gt=\nt_-\oplus\he\oplus\nt_+$.

Let $\hat{\Pi}=\{\alpha_0,\alpha_1,\dots,\alpha_\ell\}\subset
(\he)^*$ be the set of simple roots of $\gt$. Usual extensions of bilinear
forms $\langle\cdot,\cdot\rangle$ onto $\he$ and $(\he)^*$ are
denoted by the same symbols (we take $\langle c,d \rangle=1$).
Define fundamental weights $\Lambda_i\in (\he)^*$ by $\langle
\Lambda_i,\alpha_j \rangle=\delta_{ij}$ and $\Lambda_i(d)=0$,
$i,j=0,\dots,\ell$.

Let $V$ be a highest weight module for affine Lie algebra $\gt$.
Then $V$ is generated by a highest weight vector $v_\Lambda$ such
that
\begin{eqnarray*}
h\cdot v_\Lambda & = & \Lambda(h) v_\Lambda,\quad\textrm{for}\ h\in\he, \\
x\cdot v_\Lambda & = & 0,\quad\textrm{for}\ x\in\nt_+,
\end{eqnarray*}
for $\Lambda\in (\he)^*$. Module $V$ is a direct sum of weight
subspaces $V_\mu=\{v\in V \,|\,h\cdot V =\mu(h) v \textrm{ for }
h\in\he\},\,\mu\in\he$.

Standard (i.e. integrable highest weight) $\gt$-module $L(\Lambda)$
is an irreducible highest weight module, with the highest weight
$\Lambda$ being dominant integral, i.e.
$$\Lambda=k_0 \Lambda_0+k_1 \Lambda_1+\dots+k_\ell \Lambda_\ell,$$
where $k_i\in\Z_+$, $i=0,\dots,\ell$. The central element $c$ acts on
$L(\Lambda)$ as multiplication by scalar
$$k=\Lambda(c)=k_0+k_1+\dots+k_\ell,$$
which is called the level of the module $L(\Lambda)$.

\section{Feigin-Stoyanovsky's type subspace}

Vector $v\in\h$ is said to be \emph{cominuscule} if
$$\{\alpha(v)\,|\,\alpha \in R\}=\{-1,0,1\}.$$
Similarly, weight $\omega\in P$ is said to be \emph{minuscule} if
$$\{ \langle\omega,\alpha\rangle \,|\,\alpha \in R\}=\{-1,0,1\}.$$
A dominant integral weight $\omega\in P^+$
is minuscule if and only if
$$\langle\omega,\theta\rangle = 1.$$
There is a $1-1$ correspondence between cominuscule vectors and minuscule weights, and there exists finitely many minuscule weights. More precisely,
a vector $v\in\h$ is cominuscule if and only if it is dual to some
minuscule fundamental weight $\omega$, in the sense that
$$v=\nu^{-1}(\omega),$$
for some choice of positive roots.

Fix a cominuscule vector $v\in\h$. For $\g=\gsl(\ell+1,\C)$, all
 fundamental weights are minuscule. So we can assume that the
cominuscule vector $v$ is dual to a fundamental weight
$$\omega=\omega_m,$$ for some $m \in \{1,\dots, \ell\}$.
Set
$$\Gamma =
\{\,\alpha \in R \mid \alpha(v) = 1\}=\{\,\alpha \in R \mid \langle\omega, \alpha\rangle = 1\}.$$ Then we have the induced
$\mathbb Z$-gradation of ${\mathfrak g}$:
\begin{equation}
\label{GDecomp_jed} \mathfrak g  =
\mathfrak g_{-1} \oplus \mathfrak g_0 \oplus \mathfrak g_1, \end{equation}
where
\begin{eqnarray*} {\mathfrak g}_0 & = & {\mathfrak h} \oplus
\sum_{\alpha(v)=0}\, {\mathfrak g}_\alpha
\\
 \displaystyle {\mathfrak g}_{\pm1} & = &
\sum_{\alpha \in \pm \Gamma}\, {\mathfrak g}_\alpha.
\end{eqnarray*}
Subalgebras ${\mathfrak g}_1$ and ${\mathfrak g}_{-1}$ are
commutative, and ${\mathfrak g}_0$ acts on them by adjoint action.
The subalgebra ${\mathfrak g}_0$ is reductive with semisimple part
${\mathfrak l}_0=[{\mathfrak g}_0,{\mathfrak g}_0]$ of the type
$A_{m-1}\times A_{\ell-m}$; as a root basis one can take
$\{\alpha_1,\dots,\alpha_{m-1}\}\cup\{\alpha_{m+1},\dots,\alpha_\ell\}$,
and the center is equal to $\C v$.

Basis of the subalgebra ${\mathfrak g}_1$ can be identified with the set
of roots $\Gamma$. We will call elements  $\gamma\in\Gamma$ {\em
colors} and the set $\Gamma$ {\em the  set of colors}.
For $\omega=\omega_m$, the set of colors is
$$\Gamma=\{\gamma_{ij}\,|\, i=1,\ldots,m; j=m,\ldots,\ell\}$$
where
\begin{equation}\label{GammaIJ_jed} \gamma_{ij}=\alpha_i+\cdots+\alpha_m+\cdots+\alpha_j.\end{equation}
The maximal root $\theta$ is equal to $\gamma_{1\ell}$.

We picture the
set of colors $\Gamma$ as a rectangle with
row-indices $1,\dots,m$ and column-indices $m,\dots,\ell$, like in the figure \ref{Gamma_fig}.

Similarly, one also has the induced $\Z$-gradation of  affine Lie
algebra $\gt$:
\begin{eqnarray*}
\gt_0 & = & {\mathfrak g}_0\otimes\C [t,t^{-1}]\oplus \C c \oplus \C d,\\
\gt_{\pm 1} & = & {\mathfrak g}_{\pm 1}\otimes\C [t,t^{-1}],\\
\gt & = & \gt_{-1} + \gt_0 + \gt_1.
\end{eqnarray*}
As above, $\gt_{-1}$ and $\gt_1$ are commutative subalgebras, and
$\gt_1$ is a $\gt_0$-module.

For a dominant integral weight $\Lambda$, we define a
\emph{Feigin-Stoyanovsky's type subspace}
$$W(\Lambda)=U(\gt_1)\cdot v_\Lambda\subset L(\Lambda).$$

Our objective is to find a combinatorial basis of  $W(\Lambda)$. Set
$$\gt_1^+=\gt_1\cap \nt_+,\, \gt_1^-=\gt_1\cap \nt_-.$$ Then we have
$$W(\Lambda)=U(\gt_1^-)\cdot v_\Lambda.$$
By Poincar\'e-Birkhoff-Witt theorem, we have a spanning set of
$W(\Lambda)$ consisting of monomial vectors
\begin{equation}\label{PBWgen_jed}\{x_{\gamma_t}(-r_t)\cdots
x_{\gamma_2}(-r_2) x_{\gamma_1}(-r_1) v_\Lambda\,|\, t\in\Z_+;
\gamma_j\in \Gamma, r_j\in\N\}.\end{equation}

In the end, we'll say a few words about notation.  Elements of the
spanning set \eqref{PBWgen_jed} can be identified with monomials from
$U(\gt_1)=S(\gt_1)$. Because of this we often refer to
elements of $\{x_\gamma(-r) \mid \gamma\in\Gamma,r\in \Z\}$ in
$\gt_1$ as to {\em variables}, {\em elements} or {\em factors}
of a monomial.

Monomials from $S(\gt_1)$ can be identified with {\em colored
partitions}. Let
$\pi:\{x_\gamma(-r) \mid \gamma\in\Gamma,r\in \Z\}\to \Z_+$ be a
colored partition (cf. [P1], section 3). The corresponding monomial
$x(\pi)\in S(\gt_1)$ is
$$x(\pi)=x_{\gamma_t}(-r_t)^{\pi(x_{\gamma_t}(-r_t))}\cdots x_{\gamma_1}(-r_1)^{\pi(x_{\gamma_1}(-r_1))}.$$
From this identification we take notation $x(\pi)$ for the
monomials from $S(\gt_1)$. It will be convenient to define some new monomials
by using this identification. Also, our combinatorial conditions for the basis elements will be written in terms of exponents $\pi(x_\gamma(-r))$,
which gives a parametrization of the basis by a certain generalization of the notion of $(k,\ell +1)$-admissible configurations from [FJLMM].

\section{Order on the set of monomials}

\label{uredjaj_sect}

We introduce a linear order on the set of monomials.

On the weight and root lattice, we have an order $\prec$ defined
in the standard way: for $\mu,\nu\in P$ set $\mu\prec\nu$ if $\mu-\nu$
is an integral linear combination of simple roots $\alpha_i, i=1,\dots,\ell$, with non-negative coefficients.

Next, we define a linear order $<$ on the set of colors $\Gamma$ which is an extension of the order $\prec$.
For elements of $\Gamma$, $\gamma_{i'j'}\prec\gamma_{ij}$ is equivalent to saying that
$i'\geq i$ and $j'\leq j$. The order $<$ on $\Gamma$ is defined in the following way:
$$\gamma_{i'j'}<\gamma_{ij}\textrm{\quad if \quad}\left\{\begin{array}{l}i'>i \\ i'=i,\ j'<j.  \end{array}\right.$$
It is clear that this is a linear order on the set of colors.

On the set of variables
$\{x_\gamma(-r)\,|\,\gamma\in\Gamma,\,r\in\Z\}\subset \gt_1$ we
define a linear order $<$ so that we compare degrees first, and then colors of
variables:
$$x_\alpha(-r)<x_\beta(-r') \textrm{\quad if \quad}
\left\{\begin{array}{l}
-r<-r',\\
r=r' \textrm{\quad and \quad} \alpha<\beta.
\end{array}\right. $$

Since the algebra $\gt_1$ is commutative, we can assume that the
variables in monomials from $S(\gt_1)$  are sorted ascendingly
from left to right. The order $<$ on the set of monomials is
defined as a lexicographic order, where we compare variables from
right to left (from the greatest to the lowest one). If $x(\pi)$
and $x(\pi')$ are two monomials,
\begin{eqnarray*}
x(\pi) & = & x_{\gamma_t}(-r_t) x_{\gamma_{t-1}}(-t_{t-1})\cdots
x_{\gamma_2}(-r_2)x_{\gamma_1}(-r_1),\\
x(\pi') & = & x_{\gamma_s'}(-r_s')
x_{\gamma_{s-1}'}(-r_{s-1}')\cdots
x_{\gamma_2'}(-r_2')x_{\gamma_1'}(-r_1'),
\end{eqnarray*}
then $x(\pi)<x(\pi')$ if there exist $i_0\in \N$ so that
$x_{\gamma_i}(-r_i)=x_{\gamma_i'}(-r_i'),\ \textrm{for all}\ i <
i_0,$ and either $i_0=t+1\leq s$ or
$x_{\gamma_{i_0}}(-r_{i_0})<x_{\gamma_{i_0}'}(-r_{i_0}')$.

This monomial order is compatible with multiplication:
\begin{prop}[{[T]}] \label{uredjaj}
Let $$x(\pi_1)\leq x(\mu_1)\quad \textrm{and} \quad x(\pi_2) \leq
x(\mu_2).$$ Then $$x(\pi_1)x(\pi_2) \leq x(\mu_1)x(\mu_2),$$ and if
one of the first two inequalities is strict, then the last one is
also strict.
\end{prop}

For a monomials $x(\pi)\in S(\gt_1)$, we also define a degree
and a shape of $x(\pi)$. A \emph{degree} of a monomial is equal to the
sum of degrees of its variables. For
$$x(\pi) = x_{\gamma_t}(-r_t) x_{\gamma_{t-1}}(-r_{t-1})\cdots
x_{\gamma_2}(-r_2)x_{\gamma_1}(-r_1),$$ its degree is equal to
$-r_1-r_2-\dots-r_t$. A \emph{shape} of a monomial is gotten from its colored
partition by forgetting colors and considering only degrees of
factors. More precisely, for a monomial $x(\pi)$ and its partition
$\pi:\{x_\gamma(-r) \mid \gamma\in\Gamma,r\in \Z\}\to \Z_+$, the
corresponding shape will be
\begin{eqnarray*}
& & s_\pi:\Z\to\Z_+,\\
& & s_\pi(r)=\sum_{\gamma\in\Gamma}\pi(x_\gamma(-r)).
\end{eqnarray*}
A linear order can also be defined on the set of shapes; we'll say
that $s_\pi<s_{\pi'}$ if there exists $r_0\in\Z$ such that
$s_\pi(r)=s_{\pi'}(r)$ for  $r<r_0$ and either $s_\pi(r_0)<
s_{\pi'}(r_0)$ and $s_\pi(r')\neq 0$ for some $r'>r_0$, or
$s_\pi(r_0)>s_{\pi'}(r_0)$ and $s_\pi(r)=0$ for $r>r_0$.

In the end, for the sake of simplicity, we introduce the following
notation:
$$x_{ij}(-r)=x_{\gamma_{ij}}(-r),$$
 for $\gamma_{ij}\in\Gamma,\, r\in\N$.

\section{Vertex operator construction}

\label{voakonstr_sect} We use the vertex operator algebra construction
of basic $\gt$-modules (i.e. standard $\gt$-modules of level
$1$). We'll sketch this construction here, details can be found in [FLM], [DL] or [LL]; see also [FK], [S].

Consider tensor products
\begin{eqnarray*}
V_P & = & M(1)\otimes \C [P],\\
V_Q & = & M(1)\otimes \C [Q];
\end{eqnarray*}
where $M(1)$ is the Fock space for the Heisenberg subalgebra
$\hz=\sum_{n\in\Z \setminus\{0\}}\h\otimes t^n \oplus \C c$, and $\C
[P]$ and $\C [Q]$ are group algebras of the weight and root lattice
with bases consisting of $\{e^\lambda\,|\,\lambda\in P\}$, and
$\{e^\alpha\,|\,\alpha\in Q\}$, respectively. We identify group elements
$e^\lambda=1\otimes e^\lambda\in V_P$.

Space $V_Q$ has a natural structure of vertex operator algebra and
$V_P$ is a module for this algebra. Vertex operators are defined as
follows:
\begin{equation}Y(e^\lambda,z)=E^-(-\lambda,z)E^+(-\lambda,z)\otimes
e^\lambda z^\lambda \epsilon(\lambda,\cdot),
\label{vop_jed}\end{equation} where $e^\lambda=1\otimes e^\lambda$
is a multiplication operator, $\epsilon_\lambda=1\otimes
\epsilon(\lambda,\cdot)$ and $\epsilon(\cdot,\cdot)$ is a
$2$-cocycle (cf. [DL]), operator $z^\lambda=1\otimes z^\lambda$,
$z^\lambda \cdot e^\mu=e^\mu z^{\langle \lambda,\mu \rangle}$ and
$$E^{\pm}(\lambda,z)=\exp \left(\sum_{m\geq 1}\lambda(\pm m) \frac{z^{\mp m}}{\pm m}\right),$$
for $\lambda\in P$.

By using vertex operators, one can define the structure of
$\gt$-module on $V_P$. For $\alpha\in R,j\in\Z$ set
$$x_\alpha(z)=Y(e^\alpha,z),$$
for a properly chosen root vector $x_\alpha$. Heisenberg
subalgebra acts on the Fock space $M(1)$ and $c$ acts as identity.
In this way $V_Q$ and $V_Q e^{\omega_j},j=1,\dots,\ell$ become standard
$\gt$-modules of level $1$ with highest weight vectors
$v_0=1$ and $v_j=e^{\omega_j},j=1,\dots,\ell$,
$$L(\Lambda_0)\cong V_Q
\quad \textrm{and}\quad L(\Lambda_j)\cong V_Q
e^{\omega_j},j=1,\dots,\ell$$ and
$$V_P\cong L(\Lambda_0)\oplus L(\Lambda_1)\oplus\dots
\oplus L(\Lambda_\ell).$$

We will also be using intertwining operators $\mathcal Y$. More
precisely we'll be needing operators
$${\mathcal Y} (e^\lambda,z)=Y(e^\lambda,z)e^{i\pi\lambda}c(\cdot,\lambda),$$
for $\lambda\in P$, where $c(\cdot,\lambda)$ is a commutator map
(cf. [DL]).

 Restrictions of ${\mathcal Y} (e^\lambda,z)$ are in fact maps between standard modules of level
 $1$: if $\lambda+\omega_i \equiv
\omega_j \mod Q$, then
\begin{equation} \label{evop_jed} {\mathcal Y} (e^\lambda,z):L(\Lambda_i)\to
L(\Lambda_j)\{z\},\end{equation} where $L(\Lambda_j)\{z\}$ is a
space of formal series with coefficients in $L(\Lambda_j)$. Here, for convenience, we've set $\omega_0=0$. Also,
for a suitable choice of $\mu\in P$, the operators ${\mathcal
Y}(e^\mu,z_2)$ will commute with $\gt_1$ (cf. [T]).

Standard modules of level $k>1$ can be viewed, by the complete
reducibility, as submodules of tensor products of basic modules;
$$L(\Lambda)\subset L(\Lambda_0)^{\otimes
k_0}\otimes\dots\otimes L(\Lambda_\ell)^{\otimes k_\ell},$$ if
$\Lambda=k_0 \Lambda_0+k_1 \Lambda_1+\dots+k_\ell \Lambda_\ell$,
$k=k_0+k_1+\dots+k_\ell$. Highest weight vector of $L(\Lambda)$ is
$$v_\Lambda=v_{0}^{\otimes k_0}\otimes\dots\otimes
v_{\ell}^{\otimes k_\ell}.$$

This all can be imbedded into $V_P^{\otimes k}$. One can also
define vertex operators corresponding to elements
$$u_1\otimes\cdots\otimes u_k\in V_P^{\otimes k}$$
as tensor products of vertex operators on the appropriate tensor
factors:
$$Y(u_1\otimes\cdots\otimes u_k,z)=Y(u_1,z)\otimes\cdots\otimes
Y(u_k,z).$$ Then $V_Q^{\otimes k}=L(\Lambda_0)^{\otimes k}$
becomes vertex operator algebra, and $V_P^{\otimes k}$ with its
subspaces $L(\Lambda_0)^{\otimes k_0}\otimes\dots\otimes
 L(\Lambda_\ell)^{\otimes k_\ell}$ become modules for this algebra.

\section{Operator $e(\omega)$}

For $\lambda\in P$, $e^\lambda$ denotes multiplication operator
$1\otimes e^\lambda$ in $V_P=M(1)\otimes \C [P]$. Set
$$e(\lambda)=e^\lambda \epsilon(\cdot,\lambda),\qquad e(\lambda):V_P\to V_P,$$
Clearly, $e(\lambda)$ is a linear bijection. Its restrictions on
basic modules are bijections from one basic module
$L(\Lambda_i)$ onto another basic module $L(\Lambda_{i'})$. From the
definition of vertex operators $Y(e^\alpha,z),\alpha\in R$ one gets
the following commutation relation
$$Y(e^\alpha,z)e(\lambda)=e(\lambda)
z^{\langle\lambda,\alpha\rangle} Y(e^\alpha,z),$$ or, in terms of
components,
\begin{equation} x_\alpha(n)
e(\lambda)=e(\lambda) x_\alpha(n+\langle\lambda,\alpha\rangle),\quad
n\in\Z. \label{komut_ea_xn_jed}\end{equation}

For standard modules of level $k>1$, one defines operator
$e(\lambda)$ on the tensor product of basic modules as a tensor
product of the appropriate operators
$$e(\lambda)=e(\lambda)\otimes\cdots\otimes e(\lambda): \otimes_{s=1}^k L(\Lambda_{i_s})\to \otimes_{s=1}^k L(\Lambda_{i_s'}).$$
Operator $e(\lambda)$ is again a linear bijection, and relation
\eqref{komut_ea_xn_jed} still holds.

For $\lambda=\omega$ and $\gamma\in \Gamma$, the relation
\eqref{komut_ea_xn_jed} becomes
$$x_\gamma(n) e(\omega)=e(\omega)
x_\gamma(n+1).$$ More generally, for a monomial $x(\pi)\in
S(\gt_1)$, $$x(\pi)e(\omega)=e(\omega) x(\pi^+),$$ where
$x(\pi^+)\in S(\gt_1)$ denotes a monomial corresponding to partition
\begin{equation} \label{xpiplus_def}
\pi^+(x_\gamma(n+1))=\pi(x_\gamma(n)).
\end{equation}

\section{Case $k=1$}

Here we'll briefly recall the main results from [T] concerning the
basis of a Feigin-Stoyanovsky's type subspace of the standard
module $L(\Lambda_i)$. It was described in terms of difference and
initial conditions.

A monomial $x(\pi)$ satisfies
difference conditions for $L(\Lambda_i)$ if the following holds:
\begin{itemize}
\item
if $x(\pi)$ contains elements $x_{pq}(-r)$ and $x_{p'q'}(-r)$, and
$\gamma_{p'q'}\leq \gamma_{pq}$, then $p'>p$ and $q'<q$,

\item  if $x(\pi)$ contains elements
$x_{pq}(-r)$ and $x_{p'q'}(-r-1)$, then $p'>p$ or $q'<q$.
\end{itemize}

\begin{figure}[ht] \caption{Difference conditions - level $1$ case}
\label{DCbasic_fig}
\begin{center}\begin{picture}(200,140)(-8,-10) \thicklines
\put(0,0){\line(1,0){180}} \put(0,0){\line(0,1){120}}
\put(180,120){\line(0,-1){120}} \put(180,120){\line(-1,0){180}}
\put(-6,111){$\scriptstyle 1$} \put(-6,99){$\scriptstyle 2$}
\put(-8,3){$\scriptstyle m$} \put(3,-8){$\scriptstyle m$}
\put(174,-8){$\scriptstyle \ell$}

\linethickness{.075mm} \multiput(88,72)(4,0){23}{\line(1,0){2}}
\multiput(88,72)(0,4){12}{\line(0,1){2}}

\thinlines \put(96,76){$\scriptscriptstyle \bullet$}
\put(98,78){\line(1,1){12}} \put(108,88){$\scriptscriptstyle
\bullet$} \put(132,94){$\scriptscriptstyle \bullet$}
\put(134,96){\line(1,1){12}} \put(144,106){$\scriptscriptstyle
\bullet$} \put(110,90){\line(4,1){24}} \put(108,97){$\scriptstyle
(-r)$}

\put(60,10){$\scriptscriptstyle \circ$}
\put(84,16){$\scriptscriptstyle \circ$}
\put(108,40){$\scriptscriptstyle \circ$}
\put(120,64){$\scriptscriptstyle \circ$}
\put(62,12){\line(4,1){24}} \put(86,18){\line(1,1){24}}
\put(110,42){\line(1,2){12}} \put(96,22){$\scriptstyle (-r-1)$}
\end{picture}\end{center}
\end{figure}
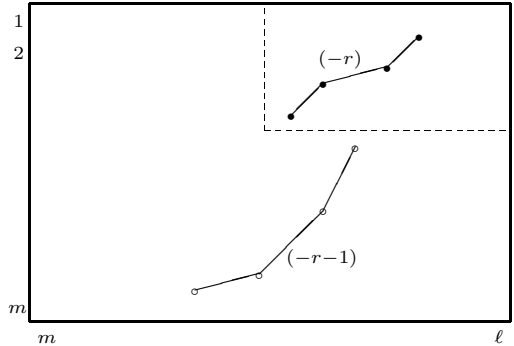

From this we can conclude that colors of the elements of the same
degree  $-r$ inside $x(\pi)$ make a descending sequence as
pictured below; appropriate row-indices strictly increase, while
column-indices strictly decrease. Colors of elements of
degree  $-r-1$ also form a decreasing sequence which is placed
below or on the left of the minimal color of elements of degree
$-r$ (see figure \ref{DCbasic_fig}).

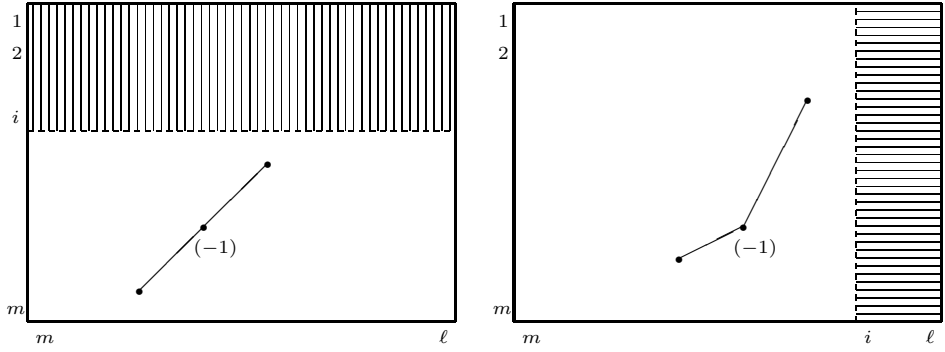
\begin{figure}[ht] \caption{Initial conditions - level $1$ case}
\label{ICbasic_fig}
\begin{center}
\begin{picture}(365,140)(-8,-10) \thicklines
\put(0,0){\line(1,0){160}} \put(0,0){\line(0,1){120}}
\put(160,120){\line(0,-1){120}} \put(160,120){\line(-1,0){160}}
\put(-6,111){$\scriptstyle 1$} \put(-6,99){$\scriptstyle 2$}
\put(-8,3){$\scriptstyle m$} \put(3,-8){$\scriptstyle m$}
\put(154,-8){$\scriptstyle \ell$}

\put(182,0){\line(1,0){160}} \put(182,0){\line(0,1){120}}
\put(342,120){\line(0,-1){120}} \put(342,120){\line(-1,0){160}}
\put(176,111){$\scriptstyle 1$} \put(176,99){$\scriptstyle 2$}
\put(174,3){$\scriptstyle m$} \put(185,-8){$\scriptstyle m$}
\put(336,-8){$\scriptstyle \ell$}

\linethickness{.075mm} \multiput(0,72)(4,0){40}{\line(1,0){2}}
\multiput(310,0)(0,4){30}{\line(0,1){2}} \put(-6,75){$\scriptstyle
i$} \put(313,-8){$\scriptstyle i$}

\multiput(2,72)(3,0){53}{\line(0,1){48}}
\multiput(310,3)(0,3){39}{\line(1,0){32}}

\thinlines \put(40,10){$\scriptscriptstyle \bullet$}
\put(42,12){\line(1,1){24}} \put(64,34){$\scriptscriptstyle
\bullet$} \put(66,36){\line(1,1){24}}
\put(88,58){$\scriptscriptstyle \bullet$}
 \put(62,26){$\scriptstyle (-1)$}

\thinlines \put(242,22){$\scriptscriptstyle \bullet$}
\put(244,24){\line(2,1){24}} \put(266,34){$\scriptscriptstyle
\bullet$} \put(268,36){\line(1,2){24}}
\put(290,82){$\scriptscriptstyle \bullet$}
 \put(264,26){$\scriptstyle (-1)$}

\end{picture}\end{center}
\end{figure}

A monomial $x(\pi)$ satisfies
initial conditions for $L(\Lambda_i)$ if it doesn't contain any element $x_{pq}(-1)\in
\gt_1$ such that $x_{pq}(-1)v_i=0$. In the case $0\leq i\leq m$, $x_{pq}(-1)v_i=0$ if $p\leq i$; in the case $m \leq i \leq \ell$,
$x_{pq}(-1)v_i=0$ if $q\geq i$. Hence, initial conditions imply that the sequence of colors of elements of
degree $-1$ lies below the $i$-th row (if $0\leq i\leq m$), or on
the left of the $i$-th column (for $m \leq i \leq \ell$); see figure \ref{ICbasic_fig}.

\section{Difference and initial conditions}

\subsection{Relations}

As in the case $k=1$ (cf. [T]) we will first find relations between
fields $x_\gamma(z), \gamma\in\Gamma$ on the standard module
$L(\Lambda)$. By equating coefficients of powers of $z$ in these
relations, we'll get relations between monomials. We'll then
identify the minimal monomial among these, the so called {\it
leading term} of a relation, and exclude from the spanning set all
monomials that contain leading terms. Difference conditions
combinatorially describe such monomials in terms of exponentials of
factors $x_{rs}(\-j)\in \gt_1^-$.

To obtain relations, we start from a consequence of
Frenkel-Kac-Segal formula for $L(\Lambda)$:
$$x_\theta(z)^{k+1}=0.$$
Since $\gt_1$ is commutative, the product of fields on the left side of
is a vertex-operator corresponding to the element $x_\theta(-1)^{k+1}1$
in $L(k\Lambda_0)\subset L(\Lambda_0)^{\otimes k}$, and the relation above is equivalent to
 $$x_\theta(-1)^{k+1}1=0.$$

By acting on this relation with $y\in \lgg_0$, one gets 
\begin{eqnarray} 0 & = & y\cdot
x_\theta(-1)^{k+1}1\,=\,[y,x_\theta(-1)^{k+1}]1 +
x_\theta(-1)^{k+1}y\cdot 1 \label{erel5_jed}\\
& = &[y,x_\theta(-1)^{k+1}]1\nonumber.
\end{eqnarray}
 A commutator $[y,x_\theta(-1)^{k+1}]$ is again an element of
 $S(\g_1(-1))\subset S(\gt_1)$ and hence from the above equality
 we obtain another relation between fields $x_\gamma(z),\,
 \gamma\in\Gamma$.
Similarly, if we start with any relation between fields
$x_\gamma(z),\,\gamma\in\Gamma$ on $L(\Lambda)$, by adjoint action
of $\lgg_0$ we would obtain a new relation between these fields.

Thus we have to study a subrepresentation $V\subset
S(\gt_1(-1))$ of $\lgg_0$ generated by a singular vector
$x_\theta(-1)^{k+1}$,
\begin{equation} \label{erel1_jed}
V=U(\lgg_0)\cdot x_\theta(-1)^{k+1}\subset S(\gt_1(-1))\subset U(\gt_1).
\end{equation}
Action of $\lgg_0$ is defined by the adjoint action of $\lgg_0$ on
$\gt_1$ (i.e. on $\g_1$).

Generally, the algebra $\lgg_0$ is a direct sum of two simple
subalgebras,
$$\lgg_0=\lgg_0'\oplus \lgg_0'',$$
with the first one being of type $A_{m-1}$, and the second one of
type $A_{\ell-m}$. The Chevalley basis of $\g$ can be chosen such
that the adjoint action of $\lgg_0$ is given by:
\begin{eqnarray}
\label{erel2_jed} [x_{-\alpha_i},x_{\gamma_{pq}}] & = &
\delta_{ip}x_{\gamma_{p+1,q}},\\ \,[x_{\alpha_i},x_{\gamma_{pq}}] & =
& \delta_{i,p-1}x_{\gamma_{p-1,q}},\nonumber
\end{eqnarray}
if $i=1,\dots,m-1$,
\begin{eqnarray}
\label{erel3_jed} [x_{-\alpha_i},x_{\gamma_{pq}}] & = &
\delta_{iq}x_{\gamma_{p,q-1}},\\ \,[x_{\alpha_i},x_{\gamma_{pq}}] &
= & \delta_{i,q+1}x_{\gamma_{p,q+1}},\nonumber
\end{eqnarray}
if $i=m+1,\dots,\ell$. One can say that the first subalgebra
$\lgg_0'$ acts by changing the row-index, and the second subalgebra
$\lgg_0''$ by changing the column-index of elements
$x_\gamma,\gamma\in\Gamma.$ 

By a simple calculation, one sees that the vector $x_\theta(-1)^{k+1}$
is the highest weight vector in $V$ of weight
$(k+1)\theta=(k+1)(\omega_1+\omega_\ell)$ for $\lgg_0$. The highest
weight representation of $\lgg_0$ can be obtained in another way, by
taking tensor products of highest weight representations its simple
subalgebras. Let $V_1$ be a highest weight representation of
$\lgg_0'$ of highest weight $(k+1)\omega_1$. It can be realized as a
subspace of homogenous polynomials of degree $k+1$ in $m$ variables,
$$V_1=S^{k+1}(x_1,\dots,x_m)\subset S(x_1,\dots,x_m).$$
The action of $\lgg_0'$ is given on the generators
$\{x_{-\alpha_i},\, x_{\alpha_i}\,|\, i=1,\dots,m-1\}$ by
\begin{eqnarray*}
x_{-\alpha_i} & \mapsto & x_{i+1}\frac{\partial}{\partial x_i},\\
x_{\alpha_i} & \mapsto & x_{i}\frac{\partial}{\partial x_{i+1}}.
\end{eqnarray*}
Similarly, let $V_2$ a highest weight representation of $\lgg_0''$
of highest weight $(k+1)\omega_\ell$. It can be realized as
homogenous polynomials of degree $k+1$ in $\ell-m+1$ variables
$$V_2=S^{k+1}(x_m,\dots,x_\ell)\subset S(x_m,\dots,x_\ell),$$ where
the action is defined on the generators $\{x_{-\alpha_i},\,
x_{\alpha_i}\,|\, i=m+1,\dots,\ell\}$ by
\begin{eqnarray*}
x_{-\alpha_i} & \mapsto & x_{i-1}\frac{\partial}{\partial x_i},\\
x_{\alpha_i} & \mapsto & x_{i}\frac{\partial}{\partial x_{i-1}}.
\end{eqnarray*}
Then
$$V\cong V_1\otimes V_2.$$

A highest weight vector in $V_1$ is $x_1^{k+1}$. From the character
formula (cf. [H]) one sees that monomials
$$x_{i_1}x_{i_2}\cdots x_{i_{k+1}},\quad 1\leq i_1 \leq
i_2\leq\dots\leq i_{k+1}\leq m,$$ constitute a basis of $V_1$ made of
weight vectors. Similarly, a monomial
$x_\ell^{k+1}$ is a highest weight vector of $V_2$, and the basis
constitutes of monomials
$$x_{j_1}x_{j_2}\cdots x_{j_{k+1}},\quad \ell\geq j_1 \geq
j_2\geq\dots\geq j_{k+1}\geq m.$$ Hence, the basis of $V_1\otimes
V_2$ is constituted by tensor products
\begin{eqnarray*}
x_{i_1} x_{i_2} \cdots x_{i_{k+1}}\otimes x_{j_1} x_{j_2} \cdots x_{j_{k+1}}, & &\\
 & & \hspace{-50ex}
1\leq
i_1 \leq i_2\leq\dots\leq i_{k+1}\leq m, \\
& & \hspace{-50ex} \ell\geq j_1 \geq j_2\geq\dots\geq j_{k+1}\geq m.
\end{eqnarray*}
The goal is to determine the corresponding basis of $V\subset
S(\gt_1(-1))$ which will then give us relations on
$L(\Lambda)$. We'll show that to the vector above corresponds a
linear combinations of all possible products of k+1 elements
$x_{pq}(-1)$ such that the (multi)set of row-indexes is equal to
$\{i_1, i_2, \dots, i_{k+1}\}$, and the (multi)set of column-indexes
is equal to $\{j_1, j_2, \dots, j_{k+1} \}$,
\begin{equation} \label{erel4_jed}
\sum_{\substack{\{p_1,\dots,p_{k+1}\}=\{i_1,\dots,i_{k+1}\}
\\\{q_1,\dots,q_{k+1}\}=\{j_1,\dots,j_{k+1}\} }} \hspace{-3ex}
C_{pq} x_{p_1 q_1}(-1) x_{p_2 q_2}(-1) \cdots x_{p_{k+1}
q_{k+1}}(-1),\end{equation} and the coefficients $C_{pq}$ are
positive integers.

One can first act on the highest weight vectors by the
$x_{-\alpha_i}$'s, $i\in\{m+1,\dots,\ell\}$, that would send
$x_1^{k+1}\otimes x_\ell^{k+1}\in V_1\otimes V_2$ to
$$x_1^{k+1}\otimes x_{j_1} x_{j_2} \cdots x_{j_{k+1}}.$$ Acting with the same $x_{-\alpha_i}$'s on $x_{1\ell}^{k+1}(-1)\in V$
will give
\begin{equation}\label{rel1_eq}
x_{1 j_1}(-1) x_{1 j_2}(-1) \cdots x_{1 j_{k+1}}(-1),
\end{equation}
 or, more precisely, some multiple of that monomial with positive
integer coefficient.  Next, one would act
on these vectors with $x_{-\alpha_{i_{k+1}-1}}
x_{-\alpha_{i_{k+1}-2}} \cdots x_{-\alpha_1}$; on one side, one gets
$$x_1^{k} x_{i_{k+1}}\otimes
x_{j_1} x_{j_2} \cdots x_{j_{k+1}}\,\in \, V_1\otimes V_2,$$ and on
the other $$\sum_{r=1}^{k+1} x_{1 j_1}(-1) \cdots x_{i_{k+1} j_r}(-1)
\cdots x_{1 j_{k+1}}(-1)\,\in\, V.$$ In the second vector for every
occurrence of index $1$ at the first place in \eqref{rel1_eq}, we have a monomial where this index was
changed to $i_{k+1}$. Monomials
$$x_{1 j_1}(-1) \cdots x_{1 j_{r-1}}(-1) x_{1 j_{r+1}}(-1)\cdots x_{1 j_{k+1}}(-1)$$
correspond to vectors
$$x_1^k \otimes x_{j_1} \cdots x_{j_{r-1}}x_{j_{r+1}} \cdots
x_{j_{k+1}}$$ from the appropriate module $V_1'\otimes V_2'$ of
the highest weight $k \theta$. Next, one would act on the vector and monomials
above by $x_{-\alpha_{i_{k}-1}} x_{-\alpha_{i_{k}-2}} \cdots
x_{-\alpha_1}$ that would change one occurrence of index $1$ into
$i_k$. Since $x_{-\alpha_s}\cdot x_{i_{k+1} j}=0$ for $s<i_{k+1}$,
the proof follows by induction on $k$.

Like in \eqref{erel5_jed}, from \eqref{erel4_jed} we obtain
\begin{equation} 
\sum_{\substack{\{p_1,\dots,p_{k+1}\}=\{i_1,\dots,i_{k+1}\}
\\\{q_1,\dots,q_{k+1}\}=\{j_1,\dots,j_{k+1}\} }} \hspace{-3ex}
C_{pq} x_{p_1 q_1}(-1) x_{p_2 q_2}(-1) \cdots x_{p_{k+1}
q_{k+1}}(-1)1=0\end{equation}
in $L(k\Lambda_0)\subset L(\Lambda_0)^{\otimes k}.$

From this we obtain the following family of relations between corresponding vertex-operators on $L(\Lambda)$:\\ for  $1\leq
i_1 \leq i_2\leq\dots\leq i_{k+1}\leq m$\quad and\quad $\ell\geq j_1 \geq j_2\geq\dots\geq j_{k+1}\geq m$,
\begin{equation}
\label{Rel_jed}
\sum_{\substack{\{p_1,\dots,p_{k+1}\}=\{i_1,\dots,i_{k+1}\}
\\\{q_1,\dots,q_{k+1}\}=\{j_1,\dots,j_{k+1}\} }} \hspace{-3ex}
C_{pq} x_{p_1 q_1}(z)x_{p_2 q_2}(z)\cdots x_{p_{k+1}
q_{k+1}}(z)=0,
\end{equation}
 where the coefficients $C_{pq}$ are
positive integers.

\subsection{Leading terms}

Fix one choice
\begin{eqnarray*}
& & \hspace{-2cm} 1\leq
i_1 \leq i_2\leq\dots\leq i_{k+1}\leq m, \\
& & \hspace{-2cm} \ell\geq j_1 \geq j_2\geq\dots\geq j_{k+1}\geq m
\end{eqnarray*}
and observe the corresponding relation \eqref{Rel_jed}. For every
$n\geq k+1$, coefficients of powers $z^{n-k-1}$ are infinite sums of
monomials: \begin{equation} \sum_{\substack{n_1+\dots+n_{k+1}=n\\
p,q}} \hspace{-2ex}C_{pq} x_{p_1 q_1}(-n_1)x_{p_2 q_2}(-n_2)\cdots
x_{p_{k+1} q_{k+1}}(-n_{k+1})=0. \label{Rel2_jed}\end{equation} In
each such sum, we identify the minimal monomial in the lexicographical
order defined in section \ref{uredjaj_sect}. We call this monomial \emph{the leading
term} of the relation. Because of the minimality, every monomial
that contains a leading term can be excluded from the spanning set
(cf. Proposition \ref{URred_prop}).

All monomials that appear in \eqref{Rel2_jed} are of the same length
$k+1$ and of the same total degree $-n$. Hence we can consider only
those monomials that are of {\em the minimal shape}, i.e. the ones in which
degrees of factors differ for at most $1$. The others will
be greater then these.

Consider first the case $n=r(k+1)$, for some $r\in\N$. In this case,
monomials of minimal shape will have all $k+1$ factors of the same
degree $-r$. So we need to find the minimal possible monomial
$$x_{p_1 q_1}(-r)x_{p_2 q_2}(-r)\cdots x_{p_{k+1}
q_{k+1}}(-r)$$
such that
\begin{eqnarray*}
\{p_1,\dots,p_{k+1}\} & = & \{i_1,\dots,i_{k+1}\},\\
\{q_1,\dots,q_{k+1}\}& = & \{j_1,\dots,j_{k+1}\}.\end{eqnarray*} Since
all factors are of the same degree, the minimal monomial will be
the one that has the minimal configuration of colors of its factors.
If we assume that factors of a monomial are sorted ascendingly
from left to right, this means that we have to choose the smallest
possible color $\gamma_{p_{k+1} q_{k+1}}$ (the greatest color in
a monomial), next the smallest possible color $\gamma_{p_{k} q_{k}}$,
and so on.

Since the row and column-indexes of colors of monomial are fixed
($\{i_1,\dots,i_{k+1}\}$ and $\{j_1,\dots,j_{k+1}\}$, resp.) the
greatest color will lie in the $i_1$-th row, and the smallest
possible of them is $\gamma_{i_1 j_{k+1}}$. The second greatest
color lies in the $i_2$-th row, the smallest possible being
$\gamma_{i_2 j_{k}}$.\footnote{If $i_1=i_2$, then $\gamma_{i_1
j_{k+1}}<\gamma_{i_2 j_{k}}$, but this shouldn't concern us,
because in this way we'll certainly obtain the smallest possible choice
of colors from the
 $i_1$-th row.} We proceed in the same manner, and obtain a monomial
\begin{equation}
\label{LT1_jed} x_{i_{k+1} j_1}(-r) \cdots x_{i_2 j_{k}}(-r) x_{i_1
j_{k+1}}(-r).
\end{equation}

\begin{figure}[ht] \caption{Colors of leading terms}
\label{LT1_fig}
\begin{center}\begin{picture}(200,140)(-8,-10) \thicklines
\put(0,0){\line(1,0){180}} \put(0,0){\line(0,1){120}}
\put(180,120){\line(0,-1){120}} \put(180,120){\line(-1,0){180}}
\put(-6,111){$\scriptstyle 1$} \put(-6,99){$\scriptstyle 2$}
\put(-8,3){$\scriptstyle m$} \put(3,-8){$\scriptstyle m$}
\put(174,-8){$\scriptstyle \ell$}

\thinlines \put(40,76){$\scriptscriptstyle \bullet$}
\put(41,78){\vector(2,-1){24}} \put(41,78){\vector(1,0){24}}
\put(41,78){\vector(0,-1){12}} \put(30,83){$\scriptstyle \gamma_{i_t
j_{k-t+2}}$} \put(66,60){$\scriptstyle \gamma_{i_{t+1} j_{k-t+1}}$}
\end{picture}\end{center}
\end{figure}
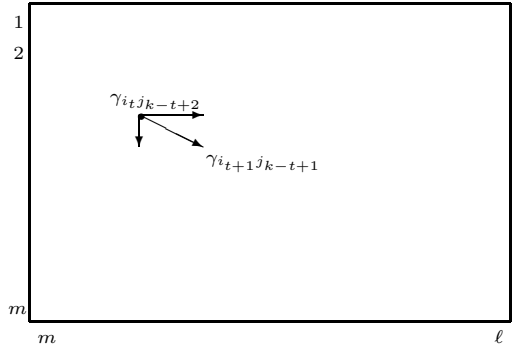

Consider now the configuration of the colors in \eqref{LT1_jed}. Each color
$\gamma_{i_{t+1} j_{k-t+1}}$ is placed on the right of, or below, or
diagonally on the right and below color $\gamma_{i_{t} j_{k-t+2}}$ (see figure \ref{LT1_fig}).

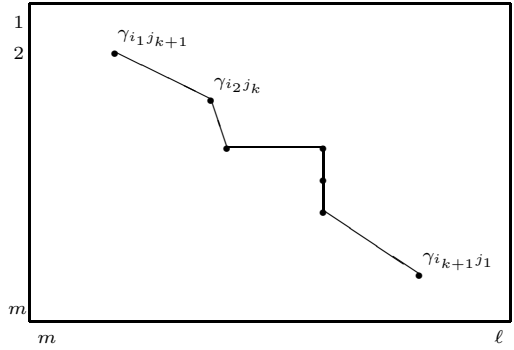
\begin{figure}[ht] \caption{Configuration of colors of leading terms}
\label{LT2_fig}
\begin{center}\begin{picture}(200,140)(-8,-10) \thicklines
\put(0,0){\line(1,0){180}} \put(0,0){\line(0,1){120}}
\put(180,120){\line(0,-1){120}} \put(180,120){\line(-1,0){180}}
\put(-6,111){$\scriptstyle 1$} \put(-6,99){$\scriptstyle 2$}
\put(-8,3){$\scriptstyle m$} \put(3,-8){$\scriptstyle m$}
\put(174,-8){$\scriptstyle \ell$}

\thinlines \put(30,100){$\scriptscriptstyle \bullet$}
\put(66,82){$\scriptscriptstyle \bullet$}
\put(72,64){$\scriptscriptstyle \bullet$}
\put(108,64){$\scriptscriptstyle \bullet$}
\put(108,52){$\scriptscriptstyle \bullet$}
\put(108,40){$\scriptscriptstyle \bullet$}
\put(144,16){$\scriptscriptstyle \bullet$}
\put(32,102){\line(2,-1){36}} \put(68,84){\line(1,-3){6}}
\put(74,66){\line(1,0){36}} \put(110,66){\line(0,-1){12}}
\put(110,54){\line(0,-1){12}} \put(110,42){\line(3,-2){36}}
\put(33,107){$\scriptstyle \gamma_{i_1 j_{k+1}}$}
\put(69,89){$\scriptstyle \gamma_{i_2 j_k}$}
\put(147,23){$\scriptstyle \gamma_{i_{k+1} j_1}$}
\end{picture}\end{center}
\end{figure}

\noindent Consequently, we conclude that colors of the leading
term lie on a diagonal path as pictured in the figure \ref{LT2_fig}

Finally, we can describe the leading terms via exponents of factors.
Let $a_{ij}^r=\pi(x_{ij}(-r))$ be an exponent of $x_{ij}(-r)$ in
$x(\pi)$. Then the leading terms correspond to solutions of
equations
\begin{eqnarray} \label{LTjed_jed}
a_{i_1 j_t}^r+ a_{i_2 j_{t-1}}^r +\dots +a_{i_t j_1}^r=k+1, & &\\
& & \hspace{-53ex}  1\leq
i_1 \leq i_2\leq\dots\leq i_t\leq m, \nonumber\\
& & \hspace{-53ex} \ell\geq j_1 \geq j_2\geq\dots\geq j_t\geq m, \nonumber \\
& & \hspace{-53ex} (i_\nu,j_\nu)\neq (i_{\nu+1},j_{\nu+1}),
\nonumber
\end{eqnarray}
where the colors $\gamma_{i_\nu j_\nu}$ lie on a diagonal path as
above.

Consider next the case $n=r(k+1)+s$. Monomials of minimal shape in this
case are built of $s$ factors of degree $-r-1$ and $k-s+1$ factors
of degree $-r$. Since factors of degree $-r$ are greater than
factors of degree $-r-1$, the leading term of \eqref{LT1_jed} can be
obtained so that first one chooses the smallest possible $(-r)$-part
of a monomial, and then the smallest possible $(-r-1)$-part. Hence
factors of degree $-r$ are placed in the last $k-s+1$ rows and the
first $k-s+1$ columns, and in the rest are placed factors of degree
$-r-1$:
\begin{eqnarray*}
(-r)\textrm{-part} & \leadsto & (i_{s+1}, i_{s+2}, \dots, i_{k+1}),\ (j_{s+1}, j_{s+2}, \dots, j_{k+1}), \\
 (-r-1)\textrm{-part} & \leadsto & (i_1, i_2, \dots, i_s),\ (j_1, j_2, \dots, j_s).
\end{eqnarray*}
We proceed as in the first case; the colors of elements of the same
degree will lie on a diagonal path as before, and the configuration of colors
is of the type pictured on the figure \ref{DCSeq_fig}.

In terms of exponents, this leading terms correspond to solutions of
equations
\begin{eqnarray*}
a_{i_1 j_s}^{r+1}+ \dots +a_{i_s j_1}^{r+1}+a_{i_{s+1} j_t}^r+ \dots
+a_{i_t j_{s+1}}^r
= k+1, & & \\
& & \hspace{-58ex} 1\leq i_1 \leq\dots\leq i_s\leq i_{s+1} \leq\dots\leq i_t\leq m, \\
& & \hspace{-58ex}  \ell\geq j_1 \geq\dots\geq j_s \geq j_{s+1}\geq\dots\geq j_t\geq m\\
& & \hspace{-58ex}  (i_\nu,j_\nu)\neq (i_{\nu+1},j_{\nu+1}),\textrm{
for }
\nu=\begin{minipage}[t]{2cm}$1,\dots,s-1,$\\$s+1,\dots,t$\end{minipage}
\end{eqnarray*}
where $a_{ij}^r=\pi(x_{ij}(-r))$, for all $r\in\N$.

In the end, observe that solutions of \eqref{LTjed_jed} are also
solutions of some equation of this type. Since the rectangles in which
colors of elements of degree $-r$ and $-r-1$ lie, intersect in
the upper-right and the lower-left corner, one can take for the
$(-r)$-part to lie in the whole $\Gamma$, and for the $(-r-1)$-part
to lie in upper-right square, on the position $1\ell$.

We'll say that a monomial  $x(\pi)\in S(\gt_1^-)$ satisfies {\em
difference conditions on} $L(\Lambda)$, or shortly, that $x(\pi)$
{\em satisfies $DC$ on} $L(\Lambda)$, if it doesn't contain a
leading term. More precisely, $x(\pi)$ satisfies difference
conditions if
\begin{eqnarray}
a_{i_1 j_s}^{r+1}+ \dots +a_{i_s j_1}^{r+1}+a_{i_{s+1} j_t}^r+ \dots
+a_{i_t j_{s+1}}^r
\leq k, & & \label{DC_eksp_jed}\\
& & \hspace{-58ex} 1\leq i_1 \leq\dots\leq i_s\leq i_{s+1} \leq\dots\leq i_t\leq m, \nonumber\\
& & \hspace{-58ex}  \ell\geq j_1 \geq\dots\geq j_s \geq j_{s+1}\geq\dots\geq j_t\geq m \nonumber\\
& & \hspace{-58ex}  (i_\nu,j_\nu)\neq (i_{\nu+1},j_{\nu+1}),\textrm{
for }
\nu=\begin{minipage}[t]{2cm}$1,\dots,s-1,$\\$s+1,\dots,t$\end{minipage} \nonumber
\end{eqnarray}
where $a_{ij}^r=\pi(x_{ij}(-r))$, for all $r\in\N$.

The following proposition was essentially proved in [P1]. The statement
follows from  the proposition \ref{uredjaj} and the minimality of leading terms. 
\begin{prop}[{[P1]}]
\label{URred_prop}
The set
\begin{equation}\{x(\pi)v_\Lambda \,|\, x(\pi)
\textrm{ satisfies DC}\}\end{equation} spans $W(\Lambda)$.
\end{prop}

\subsection{Initial conditions}

By difference conditions, in a monomial $x(\pi)$, factors of degree
$-r$, for $r>1$, are restricted by factors of degree $-r-1$ and $-r+1$.
Exceptions are factors of degree $-1$, which are restricted only
"from below", by factors of degree $-2$. Initial conditions will
play the part of restrictions "from above" on factors of degree
$-1$.

In the case $k=1$ there were no relations between monomials
consisting of factors of degree $-1$, other than those already
used for difference conditions. Thus, the initial conditions
demanded only that there are no $(-1)$-factors that annihilate
highest weight vector. However, when the level $k>1$, generally,
there are other relations between such monomials. These relations
amount to relations for difference conditions but for modules of
level lesser than $k$.

We'll say that a monomial  $x(\pi)\in S(\gt_1^-)$ satisfies {\em
initial conditions on} $L(\Lambda)$, or shortly, that $x(\pi)$ {\em
satisfies $IC$ on} $L(\Lambda)$, if
\begin{eqnarray}
& a_{i_1 j_t}^1+ a_{i_2 j_{t-1}}^1 +\dots +a_{i_t j_1}^1\leq
k_0+k_1+\dots+k_{i_t-1}+k_{j_t+1}+\dots+k_\ell, & \label{IC_jed}\\
& & \hspace{-72ex} 1\leq
i_1 \leq i_2\leq\dots\leq i_t\leq m, \nonumber\\
& &  \hspace{-72ex} \ell\geq j_1 \geq j_2\geq\dots\geq j_t\geq m, \nonumber\\
& &  \hspace{-72ex} (i_\nu,j_\nu)\neq (i_{\nu+1},j_{\nu+1}),
\nonumber
\end{eqnarray}
where $a_{ij}^1=\pi(x_{ij}(-1))$. The sum on the right side of inequality
is a sum of multiplicities of those modules of level $1$ on
which at least one $x_{i_s j_{t-s+1}}(-1)$ doesn't act as $0$.

\begin{prop}
The set
$$\{x(\pi)v_\Lambda \,|\, x(\pi) \textrm{ satisfies IC and DC}\}$$
spans $W(\Lambda)$.
\end{prop}

\begin{dokaz}
Assume that $x(\pi)$ doesn't satisfy some inequality of the type
\eqref{IC_jed} and set
$d=k_0+k_1+\dots+k_{i_t-1}+k_{j_t+1}+\dots+k_\ell+1$.  Then $x(\pi)$
contains a monomial $x(\pi')$ consisting only of  factors of degree $-1$,
such that it also doesn't satisfy that inequality. Furthermore, one
can assume that the length of $x(\pi')$ is equal to $d$. We'll show
that we can find monomials
$x(\pi_1'),\dots,x(\pi_r')$ such that $x(\pi')<x(\pi_i')$, $x(\pi_i')$'s are of the
same degree as $x(\pi')$ and $$x(\pi')v_\Lambda=\alpha_1 x(\pi_1')v_\Lambda+\dots+\alpha_t
x(\pi_t')v_\Lambda,$$ for some $\alpha_i\in\C$. Upon multiplying them with the
rest of $x(\pi)$, we'll obtain monomials $x(\pi_i)$ of the same
degree as $x(\pi)$, such that $x(\pi)<x(\pi_i)$ (cf. Proposition
\ref{uredjaj}) and
$$x(\pi)v_\Lambda=\alpha_1 x(\pi_1)v_\Lambda+\dots+\alpha_t
x(\pi_t)v_\Lambda.$$ The proof then
follows in the same way as the proof of Proposition
\ref{URred_prop} (cf. [P1,P2]).

Let $$x(\pi')=x_{i_d j_1}(-1)\cdots x_{i_1 j_d}(-1)$$ be as above.
There are $2$ possibilities:
\begin{enumerate}
\renewcommand{\labelenumi}
{(\roman{enumi})}
\item $d=k+1$.\\
In this case the initial condition is equivalent to difference
condition and the proof follows from the analogous proof for
difference conditions (cf. Proposition \ref{URred_prop}).
\item $d\leq k$.\\
Then in the tensor product $v_\Lambda=v_0^{\otimes k_0}\otimes
v_1^{\otimes k_1}\otimes \cdots \otimes v_\ell^{\otimes k_\ell}$
there is at least one $v_i$ that is annihilated by all factors of
$x(\pi')$, and $k_i>0$. Group $\Lambda_i$'s in the following way
\begin{eqnarray*}\Lambda' & = & \sum_{i=0}^{i_d-1} k_i \Lambda_i +
\sum_{i=j_d+1}^{\ell} k_i \Lambda_i,\\ \Lambda'' & =&
\sum_{i=i_t}^{j_t} k_i \Lambda_i\,=\,\Lambda - \Lambda'.
\end{eqnarray*}
Denote by $v_{\Lambda'}$ and $v_{\Lambda''}$ highest weight vectors
of  standard modules $L(\Lambda')$ and $L(\Lambda'')$. Then, by
the complete reducibility,
\begin{eqnarray*}
L(\Lambda)& \subset & L(\Lambda')\otimes L(\Lambda''),\\
v_\Lambda & = & v_{\Lambda'}\otimes v_{\Lambda''}.
\end{eqnarray*}
Since all factors of $x(\pi')$ annihilate $v_{\Lambda''}$, we have
$$x(\pi')v_\Lambda=(x(\pi')v_{\Lambda'})\otimes v_{\Lambda''}.$$
Note that $L(\Lambda')$ is a module of level $k'<k$ and $d=k'+1$.
From relations \eqref{Rel2_jed} for the module $L(\Lambda'')$ we obtain
monomials $x(\pi_1'),\dots,x(\pi_t')$ of the same length and degree,
such that $x(\pi')v_{\Lambda'}=\alpha_1
x(\pi_1')v_{\Lambda'}+\dots+\alpha_t x(\pi_t')v_{\Lambda'}$,
$\alpha_i\in\C$, and $x(\pi')<x(\pi_i')$. Also, from these relations
we see that colors of factors of monomials $x(\pi_i')$ lie in the
same rows and columns as colors of $x(\pi')$. Hence, all factors
of $x(\pi_i')$ also act as $0$ on $v_{\Lambda''}$. Consequently,
$$x(\pi')v_{\Lambda}=\alpha_1 x(\pi_1')v_{\Lambda}+\dots+\alpha_t
x(\pi_t')v_{\Lambda}.$$
\end{enumerate}
\end{dokaz}

\section{Difference and initial conditions II}

\label{UR&PU2_pogl}

We'll show how difference and initial conditions for modules of level
$k>1$ can be restated in terms of difference and initial conditions for modules of level $1$.
We are going to prove
\begin{tm}
\label{URk_tm} Let $L(\Lambda)$ be a standard module of level $k$
with highest weight vector $v_\Lambda=v_{i_1}\otimes\cdots\otimes
v_{i_k}$, where $v_{i_j}$ are highest weight vectors of
corresponding modules $L(\Lambda_{i_j})$ of level $1$. If a monomial
$x(\pi)\in S(\gt_1^-)$ satisfies difference and initial conditions
on $L(\Lambda)$, then there exists a factorization
 $$x(\pi)=x(\pi^{(1)})\cdots x(\pi^{(k)})$$
such that  $x(\pi^{(j)}) $ satisfies difference and initial
conditions on $L(\Lambda_{i_j})$.
\end{tm}

Proposition \ref{UR_ekvivalencija} will imply the converse of the
theorem. Hence, we'll have an equivalence:
\begin{kor} \label{URk_kor}
With notation as above, a monomial $x(\pi)\in S(\gt_1^-)$ satisfies
difference and initial conditions on $L(\Lambda)$ if and only if
there exists a factorization
 $$x(\pi)=x(\pi^{(1)})\cdots x(\pi^{(k)})$$
such that  $x(\pi^{(j)}) $ satisfies difference and initial
conditions on $L(\Lambda_{i_j})$.\end{kor}

\subsection{Difference conditions can be written in another way}

We first prove theorem \ref{URk_tm} in the special case $\Lambda=k
\Lambda_0$ and later for general $\Lambda$. For $\Lambda=k
\Lambda_0$, initial conditions are contained in difference
conditions, so we're only considering difference conditions on
monomials.

Define another order on the set of variables:
\begin{equation}
\label{DCuredjaj_def}
 x_{ij}(-r) \sqsubset x_{i'j'}(-r') \textrm{\quad if \quad}
\left\{\begin{array}{l} -r \leq -r'-2,\\-r=-r'-1 ;\quad i>i'\ \textrm{or}\ j<j',\\-r=-r;\quad i>i'\ \textrm{and}\ j<j'.
\end{array}\right.
\end{equation}
This is equivalent to saying that $x_{ij}(-r) < x_{i'j'}(-r')$ and
$x_{ij}(-r) x_{i'j'}(-r')$ satisfies difference conditions on
modules of level $1$. By proposition $3$ from [T], difference
conditions for modules of level 1 are conditions ``at distance 1'', i.e. relation $\sqsubset$ is
transitive, and hence it is a (strict) partial order on the set of
variables.

We look at monomials $x(\pi)\in S(\gt_1^-)$ as (multi)sets.
Then we have the following characterization of monomials satisfying difference
conditions on level $k$ modules:

\begin{prop} \label{UR_ekvivalencija}
A monomial $x(\pi)$ satisfies difference conditions on modules of
level $k$ if and only if every subset of $x(\pi)$ in which there are
no two elements comparable in the sense of $\sqsubset$, has at most
$k$ elements.
\end{prop}

\begin{dokaz}
 Let $x_{ij}(-r),x_{i'j'}(-r')\in\gt_1$ be two variables and assume $r\geq r'$. By \eqref{DCuredjaj_def}, they are incomparable if and only if
\begin{equation}
\label{DCinuredjaj_jed}
\left\{\begin{array}{l} -r = -r'; i\leq i', j\leq j',\\
-r = -r'; i\geq i', j\geq j',\\
-r=-r'-1; i\leq i', j\geq j'.
\end{array}\right.
\end{equation}
It is now clear that elements whose colors lie on a diagonal path
that was considered in \eqref{DC_eksp_jed} are mutually
incomparable in the sense of $\sqsubset$. Hence, if $x(\pi)$
doesn't satisfy difference conditions, then it has a subset of at
least $k+1$ mutually incomparable elements. Conversely, consider a
subset of $x(\pi)$ in which all elements are mutually
incomparable. By the relation \eqref{DCinuredjaj_jed}, degrees of
its elements can differ for  at most $1$. Assume that they are of
degrees $-r$ and $-r-1$. Since the elements of the same degree are
incomparable, their colors must all lie on a diagonal path like in
\eqref{LTjed_jed}. Finally, since elements of a different degree
aren't comparable these two paths are related like in
\eqref{DC_eksp_jed}.
\end{dokaz}

Notice that by proposition $3$ from [T], if $\{x_{\gamma_1}(-r_1), x_{\gamma_2}(-r_2),\dots, x_{\gamma_t}(-r_t)\}$ is a linearly ordered subset,
$$x_{\gamma_1}(-r_1) \sqsubset x_{\gamma_2}(-r_2)
\sqsubset \dots \sqsubset x_{\gamma_t}(-r_t),$$ then the monomial
$$x_{\gamma_1}(-r_1)x_{\gamma_2}(-r_2)\cdots x_{\gamma_t}(-r_t)$$
satisfies difference conditions on modules of level $1$. Thus theorem \ref{URk_tm} will be proved when we show that there exists a partition
of $x(\pi)$ into $k$
linearly ordered subsets.

\subsection{Proof of theorem \ref{URk_tm}}

Let $S$ be a finite set, $|S|=n$. Let $\sqsubset$ be a (strict)
partial order on $S$. For a subset $X\subset S$, we say that
$X$ is \emph{totaly disordered}  or \emph{discretely ordered} if elements
of $X$ are mutually incomparable, i.e. if the restriction
$\sqsubset\hspace{-.5ex}|_{X\times X}$ is an empty set.

Theorem \ref{URk_tm} now follows from the following combinatorial lemma:

\begin{lm}
Let $(S,\sqsubset)$ be a finite set with a strict partial order $\sqsubset$.
If every totaly disordered subset of $X$ consists of at most $k$ elements, then
there exists a partition of $S$ into at most $k$ linearly ordered subsets.
\end{lm}

\begin{dokaz}

 Let $l$ be the maximal cardinality of a
totaly disordered subset of $S$; $l\leq k$. We're going to show that
there is a partition of
$S$ into $l$ linearly ordered subsets.

\noindent We prove this by induction on $l$ and on the number of elements
of $S$, $n=|S|$. Distinguish 2 cases:
\begin{enumerate}
\renewcommand{\labelenumi}
{(\roman{enumi})}
\item There exists a subset $\{a_1,\dots,a_l\}\subset S$ consisting of mutually incomparable elements, such that $a_1,\dots,a_l$ aren't all maximal elements of $S$, or all minimal elements of $S$. \\
Because of the maximality of $l$, every element of $S$ is comparable to some element of $\{a_1,\dots,a_l\}$. Define subsets
$$G=\{x\in S,\ x\sqsupset a_i \ \textrm{for some}\ i\},$$
$$D=\{x\in S,\ x\sqsubset a_i \ \textrm{for some}\ i\}.$$
Since by the hypothesis $a_1,\dots,a_l$ are not all maximal elements of $S$,
nor all minimal elements of $S$, sets $G$ and $D$ are nonempty.
Then we have a partition $$S=G\cup D\cup \{a_1,\dots,a_l\}.$$  Set
$$G'=G\cup \{a_1,\dots,a_l\},$$
$$D'=D\cup \{a_1,\dots,a_l\}.$$
Sets $G'$ and $D'$ have less than $n$ elements, so by the
induction hypothesis they can be partitioned into linearly ordered
subsets. The set $\{a_1,\dots,a_l\}$ is at the same time the set
of minimal elements of $G'$, and the set of maximal elements of
$D'$. Hence, linearly ordered subsets of $G'$ end with some of the
$a_1,\dots,a_l$, while  linearly ordered subsets of $D'$ start
with some of the $a_1,\dots,a_l$. By ``gluing'' appropriate pairs
together, we get a partition of $S$ into $l$ linearly ordered
subsets.

\item The only sets with $l$ mutually incomparable elements is either the set of minimal, or the set of maximal elements of $S$. In this case we cannot
construct a partition as we did earlier because either $G$ or $D$ would be
empty. Consider 2 cases:
\begin{enumerate}
\item Assume that the only totaly disordered subset with $l$ elements
is the set of maximal elements of $S$ (analogously for minimal
elements). Denote them by $a_1,\dots,a_l$. Choose a linearly
ordered subset $\{x_1, x_2, \dots,x_r\}\subset S$ that starts with
$a_1$,
$$a_1=x_1\sqsupset x_2\sqsupset \dots\sqsupset x_r.$$  The set $S\setminus\{x_1,\dots,x_r\}$ has totaly disordered subsets of
at most $l-1$ elements, so by the induction hypothesis it can be
partitioned into $l-1$ linearly ordered subsets.
Together with $ \{x_1,x_2,\dots,x_r\}$, this gives a partition of $S$ into
$l$ linearly ordered subsets.
\item Assume that $a_1,\dots,a_l$ are all maximal, and
$b_1,\dots,b_l$ all minimal elements of $S$.\\
 Choose a linearly ordered subset $\{x_1, x_2, \dots,x_r\}\subset S$
that starts with $a_1$ and ends with some of the $b$'s, $$a_1=x_1\sqsupset x_2\sqsupset \dots\sqsupset
x_r=b_t.$$ Like in the previous case, the set
$S\setminus\{x_1,\dots,x_r\}$ has totaly disordered subsets of
at most $l-1$ elements. By the induction hypothesis it can be
partitioned into $l-1$ linearly ordered subsets, which
together with $\{x_1,x_2,\dots,x_r\}$  gives a partition of $S$ into
$l$ linearly ordered subsets.
\end{enumerate}
\end{enumerate}
\end{dokaz}

\subsection{Initial conditions}

\label{PUtrik_pogl}

We now prove the theorem \ref{URk_tm} in the general case when
$\Lambda=k_0 \Lambda_0+\dots+k_\ell \Lambda_\ell$.

First, let us recall initial conditions for a level $1$ module
$L(\Lambda_i)$, $i=1,\dots,\ell$. A monomial $x(\pi)$ satisfies initial conditions
for $L(\Lambda_i)$ if colors of elements of degree $-1$
lie below the  $i$-th row  (for
$1\leq i\leq m$), or on the left of the $i$-th column (for $m\leq i\leq \ell$).
Note that these conditions can be understood as difference conditions if
we add some imaginary elements of degree $0$ to $x(\pi)$:
for $1\leq i\leq m$ add $x_{im}(0)$ to $x(\pi)$, and for $m\leq i\leq \ell$
 add $x_{mi}(0)$ to $x(\pi)$. Then $x(\pi)$ satisfies difference and initial conditions for $L(\Lambda_i)$ if and only if this new monomial
 satisfies difference conditions for $L(\Lambda_i)$.

\begin{figure}[ht] \caption{Initial conditions in terms of difference conditions} \label{0elem_fig}
\begin{center}\begin{picture}(200,140)(-8,-10) \thicklines
\put(0,0){\line(1,0){180}} \put(0,0){\line(0,1){120}}
\put(180,120){\line(0,-1){120}} \put(180,120){\line(-1,0){180}}
\put(-6,111){$\scriptstyle 1$} \put(-6,99){$\scriptstyle 2$}
\put(-8,3){$\scriptstyle m$} \put(3,-8){$\scriptstyle m$}
\put(174,-8){$\scriptstyle \ell$}

\put(-5,74){$\scriptstyle i_t$} \put(105,-8){$\scriptstyle j_t$}

\linethickness{.075mm} \multiput(103,73)(4,0){19}{\line(1,0){2}}
\put(180,73){\line(-1,0){1}}
\multiput(103,73)(0,4){12}{\line(0,1){2}}
\multiput(109,79)(-4,0){27}{\line(-1,0){2}}
\put(0,79){\line(1,0){1}}
\multiput(109,79)(0,-4){20}{\line(0,-1){2}}

 \multiput(3,2)(6,0){18}{$\scriptscriptstyle \circ$}
 \multiput(3,8)(0,6){12}{$\scriptscriptstyle \circ$}

\thinlines

\put(5,4){\line(1,0){102}} \put(5,4){\line(0,1){72}}

\put(105,115){$\scriptscriptstyle \bullet$}
\put(118,102){$\scriptscriptstyle \bullet$}
\put(136,96){$\scriptscriptstyle \bullet$}
\put(154,90){$\scriptscriptstyle \bullet$}
\put(154,78){$\scriptscriptstyle \bullet$}
\put(107,117){\line(1,-1){13}} \put(120,104){\line(3,-1){18}}
\put(138,98){\line(3,-1){18}} \put(156,92){\line(0,-1){12}}

\put(138,102){$\scriptstyle (-1)$} \put(16,33){$\scriptstyle (0)$}

\end{picture}\end{center}
\end{figure}
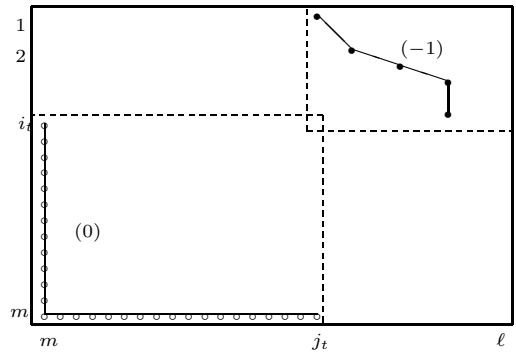

This observation generalizes to any level $k$. Let
$$\Lambda=k_0 \Lambda_0+\dots+k_\ell \Lambda_\ell,\quad k=k_0+\dots+k_\ell.$$
For every
$i=1,\dots,\ell$, we add $k_i$ elements of degree $0$ of the appropriate color to $x(\pi)$. Concretely, denote by $x(\pi')$ the monomial
$$x(\pi')=x(\pi)\cdot x_{1m}^{k_1}(0)x_{2m}^{k_2}(0)\cdots
x_{mm}^{k_m}(0)x_{m,m+1}^{k_{m+1}}(0)\cdots x_{m \ell}^{k_\ell}(0)
.$$ Colors of  $(0)$-elements of $x(\pi')$ lie on a diagonal
path as pictured on the figure \ref{0elem_fig}.
Consider difference conditions on $x(\pi')$ for
elements of degrees $-1$ and $0$. Assume that $x(\pi)$ (and $x(\pi')$)  contains elements
of degree $-1$ whose colors $\gamma_{i_1 j_t},\dots,\gamma_{i_t j_1} $ lie on a diagonal path as on the figure \ref{0elem_fig}.
Let $a_{i_1
j_t}^1,\dots,a_{i_t j_1}^1$ be the exponents of these elements. For every such choice
of $(-1)$-elements,
consider $(0)$-elements of $x(\pi')$ whose colors lie below the
$(i_t-1)$-st row and on the left of the $(j_t+1)$-st column - these are the elements
$x_{i_t m}(0),x_{i_t+1,
m}(0),\dots,x_{mm}(0),x_{m,m+1}(0),\dots,x_{mj_t}(0)$, with exponents
$k_{i_t},k_{i_t+1},\dots,k_m,\dots,k_{j_t}$, respectively. By difference conditions \eqref{DC_eksp_jed}
for $x(\pi')$, we have
$$a_{i_1 j_t}^1+\dots+a_{i_t j_1}^1+k_{i_t}+\dots+k_{j_t} \leq k.$$
Then
$$a_{i_1 j_t}^1+\dots+a_{i_t j_1}^1\leq k-k_{i_t}-\dots-k_{j_t}.$$
Hence
$$a_{i_1 j_t}^1+\dots+a_{i_t j_1}^1\leq k_0+k_1+\dots+k_{i_t-1}+k_{j_t+1}+\dots+k_\ell.$$
So, we've obtained initial conditions for $W(\Lambda)$ (cf.
\eqref{IC_jed}). We've proved
\begin{prop}
Let $x(\pi)$ and $x(\pi')$ be as above.
Then $x(\pi)$ satisfies difference and initial conditions
for $W(\Lambda)$ if and only if
$x(\pi')$ satisfies difference conditions.
\end{prop}

If $x(\pi')$ satisfies difference conditions then there exists a
partition of $x(\pi')$ into $k$ linearly ordered subsets. Elements
of degree $0$ are mutually incomparable, so they will lie in
different subsets of the partition. More precisely, they will be
the maximal elements of the corresponding subsets. By removing
these $(0)$-elements from subsets of the partition, we get a
partition of $x(\pi)$ into $k$ linearly ordered subsets. Moreover,
subsets that've contained $(0)$-elements corresponding to modules
$L(\Lambda_i), i=0,\dots,\ell$, satisfy initial conditions on
these modules. This finishes the proof of the theorem \ref{URk_tm}
in the general case.

\section{Proof of linear independence}

\label{LinNez_pogl}

\subsection{Intertwining operators}

In the level $1$ case, the main technical tool in the proof of
linear independence of the spanning set of $W(\Lambda)$ was the following
proposition (cf. [T])

\begin{prop} \label{OpIsp1_prop}
Suppose that a monomial  $x(\pi)$ satisfies difference and initial
conditions for a level $1$ standard module $L(\Lambda_i)$. Write
$x(\pi)=x(\pi_1)x(\pi_2)$, where $x(\pi_1)$ is a $(-1)$-part of a
monomial, and $x(\pi_2)$ the rest of the monomial. Then there
exists a coefficient $w(\mu)$ of an intertwining operator
${\mathcal Y}(e^\mu,z)$,
$$w(\mu):L(\Lambda_i) \to L(\Lambda_{i'})$$
for some $i'\in\{0,\dots,\ell\}$, such that:
\begin{itemize}
\item $w(\mu)$ commutes with $\gt_1$,
\item $w(\mu)x(\pi_1)v_i = C e(\omega)
v_{i'},\quad C\in\C^\times$,
\item $x(\pi_2^+)$ satisfies IC and DC for $L(\Lambda_{i'})$,
\item if  $x(\pi')$ has a $(-1)$-part $x(\pi_1')$  greater than $x(\pi_1)$, then $w(\mu)x(\pi')v_i=0$.
\end{itemize}
\end{prop}

By using theorem \ref{URk_tm}, we are able to generalize this proposition for higher level standard modules. Let $L(\Lambda)$ be a
standard module of level $k$, with the highest
weight vector $v_\Lambda=v_{i_1}\otimes\cdots\otimes v_{i_k}$.

Fix a monomial $x(\pi)$ that satisfies difference and initial conditions for $L(\Lambda)$. Let
$$x(\pi)=x(\pi_2) x(\pi_1)$$
be a factorization of $x(\pi)$ such that $x(\pi_1)$ is a $(-1)$-part, and $x(\pi_2)$ is the rest of the monomial $x(\pi)$. By theorem
\ref{URk_tm}, there exists a factorization
 $$x(\pi)=x(\pi^{(1)})\cdots x(\pi^{(k)}),$$
such that $x(\pi^{(j)})$ satisfies difference and initial conditions for $W(\Lambda_{i_j})$. Furthermore, this induces the corresponding
factorizations of $x(\pi_1)$ and $x(\pi_2)$:
\begin{eqnarray*}x(\pi_1) & = & x(\pi_1^{(1)})\cdots x(\pi_1^{(k)}),\\
x(\pi_2) & = & x(\pi_2^{(1)})\cdots x(\pi_2^{(k)}).
\end{eqnarray*}
By the proposition \ref{OpIsp1_prop}, there exist coefficients of intertwining operators $w(\mu_j),\,j=1,\dots,k$ such that
$$x(\pi_1^{(j)})v_{i_j} \xrightarrow{w(\mu_j)} C^{(j)} e(\omega) v_{i_j'},\
C^{(j)}\in\C^{\times}.$$ Let $\Lambda'=\Lambda_{i_1'}+\dots+\Lambda_{i_k'}$, and define an operator $w:L(\Lambda)\to L(\Lambda')$ with
$$w=w(\mu_1)\otimes\dots\otimes w(\mu_k).$$
Let
$$v_{\Lambda'}=v_{i_1'}\otimes\cdots\otimes v_{i_k'}$$
be the highest weight vector of $L(\Lambda')$. Then
$$x(\pi_1^{(1)})v_{i_1}\otimes\cdots\otimes x(\pi_1^{(k)})v_{i_k} \stackrel{w}{\longrightarrow} C e(\omega) v_{\Lambda'},\
C\in\C^{\times}.$$ Since by the proposition \ref{OpIsp1_prop}, $x(\pi_2^{(j)+})$ satisfy difference and initial conditions for
$L(\Lambda_{i_j'})$, then, by corollary \ref{URk_kor}, $x(\pi_2^+)$ also satisfies difference and initial conditions for $L(\Lambda')$.

Since  $x(\pi_1^{(1)})v_{i_1}\otimes\cdots\otimes
x(\pi_1^{(k)})v_{i_k}$ is only one of the summands that we get by acting
with $x(\pi_1)$ on the tensor product
$v_\Lambda=v_{i_1}\otimes\dots\otimes v_{i_k}$, we need to see what happens with other summands of $x(\pi_1)v_\Lambda$ when we
act on them with the operator $w$?
The other summands of
$x(\pi)v_\Lambda$ come from other factorizations of $x(\pi)$.
Let
$$x(\pi)=x(\nu^{(1)})\cdots x(\nu^{(k)}),$$
be another factorization of $x(\pi)$, with induced
factorizations of $x(\pi_1),x(\pi_2)$. Assume that $w$ doesn't
anihilate $x(\nu_1^{(1)})v_{i_1}\otimes\cdots\otimes
x(\nu_1^{(k)})v_{i_k}$. Since by proposition \ref{OpIsp1_prop},  $x(\pi_1^{(j)})$ is maximal for $w(\mu_j)$, we have
$$x(\nu_1^{(j)})\leq x(\pi_1^{(j)}),\quad j=1,\dots,k.$$
If some of these inequalities were strict, then by the proposition
\ref{uredjaj}, we would have $x(\pi)<x(\pi),$ which is a contradiction.
Hence, all the factors must be equal,
$$x(\nu_1^{(j)})= x(\pi_1^{(j)}),\quad j=1,\dots,k.$$
We conclude that the operator $w$ will not annihilate only those summands of
$x(\pi_1)v_\Lambda$ that are equal to
$x(\pi_1^{(1)})v_{i_1}\otimes\cdots\otimes x(\pi_1^{(k)})v_{i_k}$,
and furthermore,
$$w\cdot x(\pi_1)v_\Lambda=C\cdot( w(\mu_1)x(\pi_1^{(1)})v_{i_1}\otimes\cdots\otimes w(\mu_k)x(\pi_1^{(k)})v_{i_k}) = C' e(\omega) v_{\Lambda'},$$
for some $C,C'\in\C^{\times}.$

In a similar fashion we show that the operator $w$ annihilates all
$x(\nu)v_\Lambda$ whose $(-1)$-part $x(\nu_1)$ is greater than
$x(\pi_1)$. If $w$ doesn't annihilate $x(\nu)v_\Lambda$, then
there exists a factorization
$$x(\nu)=x(\nu^{(1)})\cdots x(\nu^{(k)}),$$
such that $x(\nu^{(j)})v_{i_j}$ aren't annihilated by $w(\mu_j)$.
We also have the induced factorization of the $(-1)$-part $x(\mu_1)$.
By the proposition \ref{OpIsp1_prop}, we have
$$x(\nu_1^{(j)})\leq x(\pi_1^{(j)}),\quad j=1,\dots,k,$$
and by the proposition \ref{uredjaj}, we conclude
$$x(\nu_1)\leq x(\pi_1).$$

We have proved a generalization of the proposition
\ref{OpIsp1_prop}:
\begin{prop} \label{OpIspk_prop}
Suppose that a monomial  $x(\pi)$ satisfies difference and initial
conditions for $L(\Lambda)$. Then there
exists an operator $w:L(\Lambda) \to L(\Lambda'),$
where $L(\Lambda')$ is another standard module of the same level, such that:
\begin{itemize}
\item $w$ commutes with $\gt_1$,
\item $w\cdot x(\pi_1)v_\Lambda= C e(\omega)
v_{\Lambda'},\quad C\in\C^\times$,
\item $x(\pi_2^+)$ satisfies IC and DC for $L(\Lambda')$,
\item $x(\pi_1)$ is maximal for $w$, i.e. all monomials  $x(\pi')$ such that $w(\mu)x(\pi')v_\Lambda\neq 0$, have their $(-1)$-part $x(\pi_1')$
smaller or equal to $x(\pi_1)$.
\end{itemize}
\end{prop}

\subsection{Proof of linear independence}

Before we proceed with the proof of linear independence, we'll
introduce some more notation, and change a bit the existing one.

For a  monomial $x(\pi)$, set $x(\pi_r)$ to be the $(-r)$-part of $x(\pi)$, and $x(\pi_r)=1$ if $x(\pi)$ doesn't contain any element of
degree $-r$. Then $$x(\pi)=x(\pi_n)x(\pi_{n-1})\cdots x(\pi_1),$$
if $x(\pi)$ consists of elements of degree greater than or equal to $-n$.
Note that the order on the set of monomials is compatible with the order
on the ''homogenous parts´´: if
\begin{eqnarray*}
x(\pi) & = & x(\pi_n)x(\pi_{n-1})\cdots x(\pi_1),\\
x(\pi') & = & x(\pi_n')x(\pi_{n-1}')\cdots x(\pi_1'),
\end{eqnarray*}
then $x(\pi')<x(\pi)$ if and only if \begin{eqnarray*}
x(\pi_1') & = & x(\pi_1)\\
& \vdots &\\
x(\pi_r') & = & x(\pi_r)\\
x(\pi_{r+1}') & < & x(\pi_{r+1})
\end{eqnarray*}
for some $r$.

Denote by $x(\pi^{+r})$ a monomial obtained from $x(\pi)$ by raising a degree by $r$ for all elements of degree $-j<-r$, and
omitting elements of degree $-j\leq -r$. Instead of $x(\pi^{+1})$, we can also write $x(\pi^+)$, for short.
Note that this new definition of $x(\pi^+)$ coincides with the old one
if degrees of elements of $x(\pi)$ are less or equal to $-2$ (cf. \eqref{xpiplus_def}).

We prove linear independence by induction. Let
\begin{equation} \label{elinzav_jed} \sum c_\pi
x(\pi)v_\Lambda=0.\end{equation} Assume that all monomials in \eqref{elinzav_jed} have elements of degree greater or equal to $-n$.
Fix a monomial $x(\pi)$ in \eqref{elinzav_jed} and assume that
$$c_{\pi'}=0 \textrm{\quad for \quad} x(\pi')<x(\pi).$$
We need to prove that $c_\pi=0$.

By the proposition \ref{OpIspk_prop}, there exists an operator $w_1:L(\Lambda)\to L(\Lambda') $ such that
\begin{itemize}
\item $w_1$ commutes with $\gt_1$,
\item $w_1\cdot x(\pi_1) v_\Lambda = C_1 e(\omega) v_{\Lambda'},\quad C_1\in\C^\times$,
\item $x(\pi^+)$ satisfies DC and IC for $L(\Lambda')$,
\item $w_1\cdot x(\pi')v_\Lambda=0$\quad for \quad $x(\pi_1')>x(\pi_1)$.
\end{itemize}
By acting with the operator $w_1$ on the relation \eqref{elinzav_jed}, we get
\begin{eqnarray*}
0 & = & w_1 \sum c_{\pi'}x(\pi') v_\Lambda\\
 & = & w_1 \sum_{\pi_1'>\pi_1}c_{\pi'}x(\pi') v_\Lambda +w_1 \sum_{\pi_1'<\pi_1}c_{\pi'}x(\pi')
 v_\Lambda+ w_1 \sum_{\pi_1'=\pi_1}c_{\pi'}x(\pi') v_\Lambda
\end{eqnarray*}
The first sum is annihilated by $w_1$ because of the maximality of $x(\pi_1)$ for $w_1$, the second sum is equal to $0$ by the induction hypothesis.
We obtain
\begin{eqnarray*}
0 & = & w_1 \sum_{\pi_1'=\pi_1}c_{\pi'}x(\pi') v_\Lambda\\
  & = & \sum_{\pi_1'=\pi_1}c_{\pi'}x(\pi_n')\cdots x(\pi_2') C_1 e(\omega) v_{\Lambda'}\\
 & = & C_1 e(\omega) \sum_{\pi_1'=\pi_1}c_{\pi'}x(\pi'^+) v_{\Lambda'}
\end{eqnarray*}
Since  $e(\omega)$ is an injection, we get
\begin{equation}
\label{elinzav2_jed} \sum_{\pi_1'=\pi_1}c_{\pi'}x(\pi'^+)
v_{\Lambda'}=0.\end{equation}

Now, for $x(\pi^+)$ there exists an operator $w_2:L(\Lambda')\to L(,\Lambda'')$
such that
\begin{itemize}
\item $w_2$ commutes with $\gt_1$,
\item $w_2\cdot x(\pi_2^+) v_{\Lambda'} = C_2 e(\omega) v_{\Lambda''},\quad C_2\in\C^\times$,
\item $x(\pi^{+2})$ satisfies DC and IC for $L(\Lambda'')$,
\item $w_2\cdot x(\pi'^+)v_{\Lambda'}=0$ \quad if \quad $x(\pi_2')>x(\pi_2)$.
\end{itemize}
Upon acting with $w_2$  on the relation \eqref{elinzav2_jed}, we get
\begin{eqnarray*}
0 & = & w_2 \sum_{\pi_1'=\pi_1} c_{\pi'}x(\pi') v_{\Lambda'}\\
 & = & w_2 \sum_{\begin{subarray}{c} \pi_1'=\pi_1\\ \pi_2'>\pi_2 \end{subarray}}c_{\pi'}x(\pi') v_{\Lambda'} +w_2
 \sum_{\begin{subarray}{c} \pi_1'=\pi_1 \\ \pi_2'<\pi_2\end{subarray}}c_{\pi'}x(\pi')
 v_{\Lambda'}+ w_2 \sum_{\begin{subarray}{c} \pi_1'=\pi_1 \\ \pi_2'=\pi_2\end{subarray}}c_{\pi'}x(\pi')
 v_{\Lambda'}
\end{eqnarray*}
As before, the first two sums are equal to $0$ because of the action of $w_2$ and of the induction hypothesis. We obtain
\begin{eqnarray*}
0 & = & w_2 \sum_{\begin{subarray}{c} \pi_1'=\pi_1 \\ \pi_2'=\pi_2\end{subarray}}c_{\pi'}x(\pi')v_{\Lambda'} \\
  & = & \sum_{\begin{subarray}{c} \pi_1'=\pi_1 \\ \pi_2'=\pi_2\end{subarray}}c_{\pi'}x(\pi_n')\cdots x(\pi_2') C_2 e(\omega) v_{\Lambda''}\\
 & = & C_2 e(\omega) \sum_{\begin{subarray}{c} \pi_1'=\pi_1 \\ \pi_2'=\pi_2\end{subarray}}c_{\pi'}x(\pi'^{+2}) v_{\Lambda''}.
\end{eqnarray*}
Since $e(\omega)$ is an injection, we get
$$\sum_{\begin{subarray}{c} \pi_1'=\pi_1 \\ \pi_2'=\pi_2\end{subarray}}c_{\pi'}x(\pi'^{+2})
v_{\Lambda''}=0.$$ We proceed inductively; after $n$ steps, we obtain
$$0=\sum_{\begin{subarray}{c} \pi_1'=\pi_1 \\ \pi_2'=\pi_2 \\ \vdots \\  \pi_n'=\pi_n \end{subarray}}c_{\pi'}x(\pi'^{+n}) v_{\Lambda^{(n)}}
=c_\pi x(\pi^{+n}) v_{\Lambda^{(n)}}=c_\pi v_{\Lambda^{(n)}}$$
and we can conclude that $c_\pi=0$.

Hence we have proved
\begin{tm}\label{FSbaza_tm}
The set
$$\{x(\pi)v_\Lambda \,|\, x(\pi) \ \textrm{satisfies DC and IC for}\ L(\Lambda)\}$$
is a basis for $W(\Lambda)$.
\end{tm}

\section{Basis of a standard module}

Feigin-Stoyanovsky's type subspace $W(\Lambda)$ was implicitly introduced and studied in  [P1] and [P2] where the basis of the
whole standard module $L(\Lambda)$ was constructed from a basis of this subspace. We have used this approach in [T] to construct
a basis of standard modules of level $1$, for any possible choice of $\Z$-gradation \eqref{ZGradG_jed}. By using corollary \ref{URk_tm}
we are able to extend this proof to standard modules of higher level.

Set
$$e=\prod_{\gamma\in\Gamma}e^\gamma=e^{\sum_{\gamma\in\Gamma}\gamma}.$$
Then
$$e=e^{(\ell+1)\omega}$$
(cf. relation $(26)$ in [T]). The following proposition was proven by Primc (cf.
Theorem 8.2. in [P1] or Proposition 5.2. in [P2])

\begin{prop} \label{StModGenSkup_prop} Let $L(\Lambda)_\mu$ be a weight subspace of
$L(\Lambda)$. Then there exists an integer $n_0$ such that for any
fixed $n\leq n_0$ the set of vectors
$$e^n x_{\gamma_1}(r_1)\cdots x_{\gamma_s}(r_s)v_\Lambda \in L(\Lambda)_\mu,$$
where  $s\geq
0,\,\gamma_1,\dots,\gamma_s\in\Gamma,\,r_1,\dots,r_s\in\Z$, is a
spanning set of $L(\Lambda)_\mu$. In particular,
$$L(\Lambda)=\langle e \rangle U(\gt_1)v_\Lambda.$$
\end{prop}

We'll use our results on the basis of $W(\Lambda)$ to prove
\begin{tm}
\label{StModBase_tm}
Let $L(\Lambda)_\mu$ be a weight subspace of a standard
$\gt$-module $L(\Lambda)$. Then there exists $n_0\in\Z$ such that
for any fixed $n\leq n_0$ the set of vectors
$$e^n x(\pi)v_\Lambda\in L(\Lambda)_\mu, \quad x(\pi)\textrm{ satisfies IC and DC for } L(\Lambda),$$
is a basis of $L(\Lambda)_\mu$.
Moreover, for two choices of
$n_1,n_2\leq n_0$, the corresponding two bases are connected by a
diagonal matrix.
\end{tm}

We have proven this theorem in [T] for standard modules of level $1$.
The first part of the theorem directly follows from proposition \ref{StModGenSkup_prop} and theorem \ref{FSbaza_tm}.
For the second part of the theorem in [T] we have considered a monomial $x(\mu)\in U(\gt_1^-)$ which was defined as
the maximal monomial satisfying difference and initial conditions for $L(\Lambda)$ such that its
factors are of degree greater or equal to $-r\frac{\ell+1}{m(\ell-m+1)}$, where $r$ is equal to the smallest common
multiple of $m$ and $\ell-m+1$. For simplicity, set $f=r\frac{\ell+1}{m(\ell-m+1)}$.
We've shown that the following holds
\begin{enumerate}
\renewcommand{\labelenumi}
{(\roman{enumi})}
\item $e(\omega)^f v_\Lambda=C x(\mu)v_\Lambda$, for some $C\in
\N$,
\item $f$ divides $\ell +1$,
\item if a monomial $x(\pi)$ satisfies difference  and initial conditions for $L(\Lambda)$, then so does
a monomial $x(\pi^{-f})x(\mu)$, where $\pi^{-f}$ is a partition
defined by $$\pi^{-f}(x_\gamma(-n-f))=\pi(x_\gamma(-n)),
\ \gamma\in\Gamma, n\in\Z.$$
\end{enumerate}
Then we've had $$e(\omega)^f x(\pi)v_\Lambda=x(\pi^{-f})e(\omega)^f
v_i=C x(\pi^{-f})x(\mu)v_\Lambda.$$ Since $e^{\omega} x(\pi)v_\Lambda$ and
$e(\omega) x(\pi)v_\Lambda$ are proportional, the second part of the
theorem followed.

To prove the theorem for higher levels, it is enough to construct a monomial
$x(\mu)\in U(\gt_1^-)$ that satisfies properties (i) and (iii).
Let $L(\Lambda)$ be a standard module of level $k$, with the highest
weight vector $v_\Lambda=v_{i_1}\otimes\cdots\otimes v_{i_k}$.
For each $L(\Lambda_{i_j})$, let $x(\mu^{(j)})$ be as in the previous paragraph.
Set
$$x(\mu)=x(\mu^{(1)})\cdots x(\mu^{(k)}).$$
As in the proof of proposition \ref{OpIspk_prop}, we have that
$x(\mu)$ is a maximal monomial satisfying difference and initial conditions for $L(\Lambda)$ such that its
factors are of degree greater or equal to $-f$ and
$$ x(\mu)v_\Lambda = C\cdot x(\mu^{(1)})v_{i_1}\otimes\cdots\otimes x(\mu^{(k)})v_{i_k}. $$
Because of the property (i), we have
\begin{eqnarray*}
  x(\mu)v_\Lambda & = & C \cdot e(\omega)^f v_{i_1}\otimes\cdots\otimes e(\omega)^f v_{i_k}\\
 & = & C e(\omega)^f v_\Lambda,
\end{eqnarray*}
for some $C\in\N$. Finally, by the property (iii) and the corollary \ref{URk_kor}, it follows that if a monomial $x(\pi)$
satisfies difference  and initial conditions for $L(\Lambda)$, then so does a monomial $x(\pi^{-f})x(\mu)$.

\section{Presentation of $W(\Lambda)$}

\label{Prez_pogl}

By definition Feigin-Stoyanovsky's type subspace $W(\Lambda)$ is a $\gt_1$ submodule of $L(\Lambda)$ generated by the highest-weight vector $v_\Lambda$,
$$W(\Lambda)=U(\gt_1)\cdot v_\Lambda.$$
Since the $\gt_1$ is commutative,
$U(\gt_1)\cong S(\gt_1)$, we have
$$W(\Lambda)=\C[x_\gamma(-r)\,|\,r\in\N]\cdot v_\Lambda.$$

For $1\leq i <m$ and $m<j\leq\ell$ let $\g_{ij}\subset \g_0$
be the subalgebra generated by the elements $x_{\pm\alpha_t}$, where
either $1\leq t \leq i-1$ or $j+1\leq t \leq\ell$.

Consider a polynomial algebra $\mathcal{A}=\C[x_\gamma(-r)\,|\,r\in\N,\gamma\in\Gamma]$,
which is also a $\g_0$-module. Let $J\subset\mathcal{A}$ be the ideal generated by the following sets
$$U(\g_0)\cdot\left(\sum_{\substack{r_1,\dots,r_{k+1}\leq -1 \\ r_1+\dots+r_{k+1}=n}}x_\theta(r_1)\cdots x_\theta(r_{k+1})\right),\quad \textrm{for all\ }n\in \Z_{<0},$$
and $$ U(\g_{ij})\cdot
\left(x_\theta(-1)^{k_0+\dots+k_{i-1}+k_{j+1}+\dots+k_\ell}\right),\quad
\textrm{for all\ } \begin{array}{l} i=1,\dots,m-1;\\
j=m+1,\dots,\ell.\end{array}$$
Then we have the following presentation result:

\begin{tm}
As a vector space, $W(\Lambda)$ is isomorphic to the quotient $\mathcal{A}/J$.
\end{tm}

\begin{dokaz}
Define a mapping
\begin{eqnarray*}
\varphi_0:\mathcal{A} & \to &  W(\Lambda),\\
\varphi_0:x(\pi) & \mapsto & x(\pi)\cdot v_\Lambda.
\end{eqnarray*}
Since the ideal $J$ lies in the kernel of $\varphi_0$, we can factorize $\varphi_0$
to a quotient map
$$
\varphi: \mathcal{A}/J \to  W(\Lambda).
$$ The map $\varphi$ is clearly a surjection, since $\varphi_0$ is
a surjection. We'll show that $\varphi$ is also an injection. Consider a set
$$\mathcal{B}=\{x(\pi)\,|\,x(\pi)\ \textrm{satisfies DC and IC for}\ L(\Lambda)\}\subset \mathcal{A}/J.$$
As in the proof of the proposition \ref{URred_prop}, we see that
this set spans $\mathcal{A}/J$. Since $\varphi$
bijectively maps this set onto
$$\{x(\pi)v_\Lambda\,|\,x(\pi)\ \textrm{satisfies DC and IC for}\ L(\Lambda)\}\subset W(\Lambda),$$
which is a basis of $W(\Lambda)$, we see that $\mathcal{B}$ is also linearly
independent. Hence $\varphi$ maps
a basis of $\mathcal{A}/J$ onto a basis of $W(\Lambda)$ and therefore
 $\varphi$ is a bijection.
\end{dokaz}

This kind of presentation of $W(\Lambda)$ was used in [FJMMT] in
order to obtain fermionic formulas for the character of
$W(\Lambda)$. Also, presentation of the Feigin-Stoyanovsky's
principal subspace was used in [C1,C2,CalLM1,CalLM2,CLM2], for construction of exact sequences between
different principal subspaces from which they obtained recurrence
relations for the characters of these spaces.

\end{document}